\documentclass{amsart}

\usepackage{amssymb}

\newcommand{\mfT}{\mathfrak T}
\newcommand{\umfT}{\underline{\mathfrak T}}
\newcommand{\imfT}{{\mathfrak T}^{-1}}
\newcommand{\uimfT}{\underline{{\mathfrak T}}^{-1}}
\newcommand{\umfa}{\underline{\mathfrak a}}
\newcommand{\umfb}{\underline{\mathfrak b}}
\newcommand{\dsp}{\displaystyle}
\newcommand{\eps}{\varepsilon}
\newcommand{\ue}{u^\eps}
\newcommand{\dt}{\partial_t}

\newcommand{\R}{\mathbb R}

\newcommand{\N}{\mathbb N}
\newcommand{\cN}{{{\mathcal N}}}
\newcommand{\cD}{{{\mathcal D}}}
\newcommand{\cF}{{{\mathcal F}}}
\newcommand{\cFe}{{{\mathcal F}^\eps}}
\newcommand{\ucFe}{{\underline{{\mathcal F}}^\eps}}
\newcommand{\cGe}{{{\mathcal G}^\eps}}
\newcommand{\cGeu}{{{\mathcal G}^\eps_u}}
\newcommand{\cGeuu}{{{\mathcal G}^\eps_{uu}}}

\newcommand{\cFeu}{{{\mathcal F}^\eps_u}}

\newcommand{\cL}{{\mathcal L}}
\newcommand{\cLe}{{{\mathcal L}^\eps}}
\newcommand{\init}{_{\vert_{t=0}}}
\newcommand{\uv}{\underline{v}}
\newcommand{\uc}{\underline{c}}
\newcommand{\uV}{\underline{V}}
\newcommand{\uz}{\underline{\zeta}}
\newcommand{\uT}{\underline{T}}
\newcommand{\cT}{{\mathcal T}}
\newcommand{\cQ}{{\mathcal Q}}
\newcommand{\uinit}{\underline{u}_0^\eps}
\newcommand{\ve}{v^\eps}

\newcommand{\us}{\underline{s}}

\newcommand{\tue}{\widetilde{u}^\eps}
\newcommand{\uue}{\underline{u}^\eps}
\newcommand{\uu}{\underline{u}}
\newcommand{\wthe}{\widetilde{h}^\eps}

\newcommand{\rge}{{\eps\in(0,\eps_0)}}
\newcommand{\Rge}{{0<\eps<\eps_0}}
\newcommand{\cst}{\mbox{\textnormal{Cst }}}

\newcommand{\supT}{\sup_{t\in[0,T]}}

\newcommand{\sk}{$_{k}\,$}
\newcommand{\skp}{$_{k+1}\,$}

\newcommand{\uapp}{u_{app}^\eps}
\newcommand{\uh}{\underline{h}}
\newcommand{\oper}{\mfT}
\newcommand{\uoper}{(\uh+\mu\underline{{\mathcal T}})}

\newcommand{\dive}{\mbox{\textnormal{div }}}

\theoremstyle{theorem}
\newtheorem{theorem}{Theorem}
\newtheorem{assumption}{Assumption}
\newtheorem{proposition}{Proposition}
\newtheorem{lemma}{Lemma}
\newtheorem{corollary}{Corollary}

\theoremstyle{definition}
\newtheorem{definition}{Definition}
\theoremstyle{remark}
\newtheorem{example}{Example}
\newtheorem{remark}{Remark}

\begin{document}

\newtheoremstyle{assumptionbis}
  {\topsep}
  {\topsep}
  {\itshape}
  {}
  {\bfseries}
  {'.}
  {.5em}
  {}

\theoremstyle{assumptionbis}
\newtheorem{assumptionbis}{Assumption}
\newtheorem{lemmabis}{Lemma}
\newtheorem{theorembis}{Theorem}

\title[A Nash-Moser theorem for singular evolution equations]{A Nash-Moser theorem for singular evolution equations. Application to the Serre and Green-Naghdi equations}         
\author{Borys Alvarez-Samaniego \and David Lannes}                     
\address{Universit\'e Bordeaux I; IMB et CNRS UMR 5251\\351 Cours de la Lib\'eration\\33405 Talence Cedex, France}                       
    \thanks{This work was supported by the ACI Jeunes Chercheuses et 
Jeunes Chercheurs ``Dispersion et nonlin\'earit\'es''.}                       

\begin{abstract}
We study the well-posedness of the initial value problem for a wide class
of singular evolution equations. We prove a general well-posedness theorem
under three assumptions easy to check: the first controls the singular part
of the equation, the second the behavior of the nonlinearities, and the
third one assumes that an energy estimate can be found for the linearized
system. We allow losses of derivatives in this energy estimate and therefore
construct a solution by a Nash-Moser iterative scheme. As an application to
this general theorem, we prove the well-posedness of the Serre and
Green-Naghdi equation and discuss the problem of their validity as 
asymptotic models for the water-waves equations.
\end{abstract}
 \maketitle  
\section{Introduction}

\subsection{General setting}

We investigate in this paper the local in time 
well-posedness of singular evolution
equations of the form
\begin{equation}
	\label{intro1}
	\left\lbrace
	\begin{array}{l}
	\dsp \dt \uue +\frac{1}{\eps}\cLe(t) \uue+\cFe[t,\uue]=h^\eps\\
	\uue\init=\uinit,
	\end{array}\right.
\end{equation}
where $\eps\in (0,\eps_0)$ is a parameter, $\cLe(t)$ is a linear operator,
while $\cFe[t,\cdot]$ is nonlinear. Under appropriate assumptions, we prove
that the initial value problems (IVP) (\ref{intro1})$_\Rge$ admit a
solution on a time interval $[0,T]$, with $T>0$ \emph{independent of $\eps$}.

Such a result is known in the case of quasilinear symmetric
hyperbolic systems, and provided that the linear (and singular) part
$\frac{1}{\eps}\cLe(t)$ is, say, a constant coefficient anti-adjoint 
differential operator (see e.g. \cite{TaylorM3} for the case of
classical symmetric system, and \cite{Grenier} for an extension of 
these results).\\
In the quasilinear case for instance, an essential step is the study of
the IVP associated to the linearization 
of (\ref{intro1}) around
any reference function $u$ belonging to some functional space $X$: 
if a solution $v$ to this IVP can be found
in $X$, and if  an 
\emph{energy estimate} controls the norm of $v$ in terms of the
norm of $u$, then a solution to (\ref{intro1}) can be constructed by
a standard Picard iterative scheme.\\
Our goal here is to investigate situations where this general approach fails.
In particular, it sometimes happens that the energy estimate  associated
to the linearized problem only
controls $v$ in a space strictly larger than $X$; when such a loss of 
information occurs, the standard Picard iterative scheme cannot converge. 
It is however possible, under certain assumptions, to use the iterative
scheme developed by Nash and Moser and used for the first time to solve 
the embedding problem for Riemannian manifolds \cite{Nash}. There exists now
an extensive literature 
(e.g. \cite{Hamilton,AlinhacGerard}) showing that the technique of Nash
and Moser can be used to prove an abstract implicit function theorem.\\
The implementation of a Nash-Moser iterative scheme is however very technical,
and is only used as a last recourse to solve nonlinear evolution equations,
though some recent works show that it is a useful tool 
(e.g. \cite{Poppenberg,PoppenbergSchmittWang,Lindblad1,Lindblad2,LannesJAMS,IOP}). We develop here a
Nash-Moser theorem specific to the general class of 
IVP (\ref{intro1}), which allows
us to greatly simplify the general theory (at the cost, sometimes, of
optimality -- see also \cite{SaintRaymond} for a simplified general
Nash-Moser implicit function theorem). The interest of these simplifications
is twofold: i) we can state a general well-posedness theorem for (\ref{intro1})
under three assumptions easy to check on $\cL^\eps(t)$, $\cFe$ and the linearization
of (\ref{intro1}); ii) we can also handle the presence in the equation of 
parameters and singular terms. We also show how these results can be 
used for the justification of asymptotic systems.\\
As an illustration, we solve the Serre and Green-Naghdi equations which
are two of the most widely used models in coastal oceanography 
(\cite{GreenLawsNaghdi,GreenNaghdi,BazdenkovMorozovPogutse} and, for instance,
\cite{Wei_etal,KimBaiErtekinWebster}). We also address the problem of the
relevance of these models as asymptotic models for the exact water-waves
equations.

\subsection{Organization of the paper}

We start by giving the three assumptions of our general well-posedness 
theorem for (\ref{intro1}) in Section \ref{sectass}. Section \ref{sectNM}
is devoted to the main theorem:
it is stated in Section \ref{sectstatement} and proved in
Sections \ref{sectproof} and \ref{sectprooflemmas}.\\
In Section \ref{sectfurther}, we give some generalizations and a corollary
of the theorem. The three main assumptions are weakened in 
Section \ref{sectgen}
where we allow a more complex dependence of the energy estimate on
time derivatives. In Section \ref{sectremarks}, some useful and easy
generalizations are given: a slight weakening of the three main
assumptions (\ref{restr}), the possibility of handling other
parameters than $\eps$ (\ref{other}) and of replacing the
linearization of (\ref{intro1}) by an approximate linearization (\ref{sectappr}). Finally, a corollary is given in Section \ref{sectstab}, which
gives a stability property very useful for the justification of
asymptotics to (\ref{intro1}).\\
An application of the main theorem is given in Section \ref{sectGN} where
the Serre and Green-Naghdi equations are solved 
uniformly with respect
to the so-called shallowness parameter (Section \ref{sectWP}). The results of
Section \ref{sectGN} are then used in Section \ref{sectjustif} to
address the justification of the Serre and Green-Naghdi models as
asymptotic models for the full water-waves equations.

\subsection{Notations}
- We generically denote by $C(\lambda_1,\lambda_2,\dots)$ a constant
depending on the parameters $\lambda_1,\lambda_2,\dots$; \emph{the dependence
on the $\lambda_j$ is always assumed to be nondecreasing}.\\
- If $X_1$ and $X_2$ are two Banach spaces, we denote by 
${\mathfrak L}(X_1,X_2)$ the set of all continuous linear mappings defined on $X_1$ and
with values in $X_2$.\\
- If $X$ is a Banach space and $T>0$, then $X_T$ stands for
$C([0,T];X)$, and we denote by $\vert\cdot\vert_{X_T}$ its canonical
norm.\\
- If $X_1$ and $X_2$ are two Banach spaces and 
${\mathcal F}\in C([0,T]; C^j(X_1;X_2))$, we denote by ${\mathcal F}_u$,
${\mathcal F}_{uu}$ and ${\mathcal F}_{(j)}$ the first, second and $j$-th
order derivatives of the
mapping $u\mapsto {\mathcal F}[\cdot,u]$.\\
- If $X_1$ and $X_2$ are two Banach spaces and 
${\mathcal F}\in C^j([0,T]; C(X_1;X_2))$, we denote by ${\mathcal F}^{(j)}$
 the $j$-th
order derivative of the
mapping $t \mapsto {\mathcal F}[t,\cdot]$.\\
- We denote $\Lambda:=(1-\Delta)^{1/2}$ and $H^s(\R^d)$ ($s\in\R$)
the usual Sobolev space 
$H^s(\R^d)=\{u\in {\mathcal S}'(\R^d),\vert u\vert_{H^s}<\infty\}$, 
where $\vert u\vert_{H^s}=\vert \Lambda^s u\vert_{L^2}$.
We keep this notation if $u$ is a vector or matrix with coefficients
in $H^s(\R^d)$.\\
- We use the condensed notation
\begin{equation}
	\label{nota3}
	A_s=B_s +\left\langle C_s\right\rangle_{s> \us}
\end{equation}
to say that $A_s=B_s$ if $s\leq \us$ and $A_s=B_s+C_s$ if $s> \us$.\\
- By convention, we take $\dsp \sum_{j=1}^0=0$ and $\dsp\prod_{j=1}^0=1$.

\subsection{Main assumptions}\label{sectass}

We state here three assumptions which imply the well-posedness of 
(\ref{intro1}). The first one deals  with the linear operator $\cLe$, the
second one with the
nonlinear term $\cFe$, and the last one with the well-posedness of the
linearization of (\ref{intro1}). Throughout this article, we assume that $(X^s)_{s\in\R}$ is a Banach scale 
in the following sense:
\begin{definition}
	\label{defintro1}
	We say that a family of Banach spaces 
	$((X^s),\vert\cdot\vert_s)_{s\in\R}$ is a Banach scale if:
	\begin{itemize}
	\item For all $s\leq s'$, one has $X^{s'}\subset X^s$ and
	$\vert\cdot\vert_s\leq \vert\cdot\vert_{s'}$;
	\item There exists a family of \emph{smoothing operators} 
	${\mathcal S}_\theta$ ($\theta\geq 1$) such that 
	$$
	\forall s<s',\quad\forall u\in X^{s'},
	\quad \vert (1-{\mathcal S}_\theta)u\vert_s\leq C_{s,s'}
	\theta^{s-s'}\vert u\vert_{s'}
	$$
	and
	$$
	\forall s\leq s', \quad \forall u\in X^s,
	\quad {\mathcal S}_\theta u\in X^{s'}
	\quad\mbox{ and }\quad
	\vert {\mathcal S}_\theta u\vert_{s'}\leq 
	C_{s,s'} \theta^{s'-s}\vert u\vert_s;
	$$
	\item The norms satisfy a convexity property:
	$$
	\forall s\leq s''\leq s',\quad\forall u\in X^{s'},\qquad
	\vert u\vert_{s''}\leq C_{s,s',s''}\vert u\vert_{s}^\mu
	\vert u\vert_{s'}^{1-\mu},
	$$
	where $\mu$ is given by the relation $\mu s+(1-\mu)s'=s''$.
	\end{itemize}
\end{definition}

The assumption made on the linear operator $\cLe$ is the following:
\begin{assumption}
	\label{assintro1}
	There exist $T>0$, 
	$s_0\in \R$ and $m\geq 0$ such that:
	\begin{enumerate}
	\item \label{condass11} For all $s\geq s_0$, one has
	$\cLe\in C(\R;{\mathfrak L}(X^{s+m};X^{s}))$ and 
	 $(\cLe(\cdot))_\Rge$ is  bounded in 
	$C([0,T];{\mathfrak L}(X^{s+m};X^{s}))$;
	\item \label{condass12} One can define an evolution operator 
	$U^\eps(\cdot)\in C(\R;{\mathfrak L}(X^s,X^s))$ 
	($s\geq s_0$) as
	$$
	\forall g\in X^s,\quad U^\eps(t)g:=u^\eps(t), 
	\quad \mbox{ where } \quad 
	\dt u^\eps+\frac{1}{\eps}\cLe(t)u^\eps=0,
	\qquad u^\eps\init=g,
	$$
	and $(U^\eps(\cdot))_\Rge$ is  bounded
	in $C([-T,T];{\mathfrak L}(X^s,X^s))$.
	\end{enumerate}
\end{assumption}

We can now state our assumption on the nonlinear operator $\cFe$: 
\begin{assumption}
	\label{assintro2}
	There exist $m \geq 0$, $T>0$, and $s_0\in\R$ such that for all 
	$s\geq s_0$, $\cF\in C([0,T];C^2(X^{s+m},X^s))$ and:
	\begin{enumerate}
	\item\label{condass1} For all $u\in X^{s+m }$,
	$$
	\sup_{t\in [0,T]}\vert 
	\cFe[t,u]\vert_s\leq C(s,T,\vert u\vert_{s_0+m })\vert u\vert_{s+m };
	$$
	\item\label{condass2} For all
	$u,v\in X^{s+m}$ one has 
	\begin{eqnarray*}
	\supT\vert {\mathcal F}^{\eps}_{u}[t,u]v
	\vert_s&\leq& 
	C(s,T,\vert u\vert_{s_0+m })\big(\vert v\vert_{s+m}
	+\vert u\vert_{s+m}
	 \vert v\vert_{s_0+m}\big);
	\end{eqnarray*}
	\item\label{condass3}  For all
	$u,v_1,v_2\in X^{s+m}$ one has 
	\begin{eqnarray*}
	\supT\vert {\mathcal F}^{\eps}_{uu}[t,u](v_1,v_2)
	\vert_s&\leq& 
	C(s,T,\vert u\vert_{s_0+m })\big(\vert v_1\vert_{s+m}
	\vert v_2\vert_{s_0+m}+\vert v_1\vert_{s_0+m}
	\vert v_2\vert_{s+m}\\
	& &+\vert u\vert_{s+m}
	 \vert v_1\vert_{s_0+m} \vert v_2\vert_{s_0+m}\big).
	\end{eqnarray*}
	\end{enumerate}
\end{assumption}
\begin{remark}
	The estimates
	of the assumption are uniform with respect to $\rge$ and
	called 
	\emph{tame estimates} after Hamilton \cite{Hamilton}: 
	the dependence of the r.h.s.
	 on the norms involving the index $s$ is
	\emph{linear}.
\end{remark}

Before stating  the assumption made on the linearization of 
(\ref{intro1}), let us  define the space $X^{s}_{(j)}$ ($j\in\N$) and $F^s$ as
\begin{eqnarray}
	\label{Xs}
	X^{s}_{(j)} &:=&\bigcap_{k=0}^j  C^k([0,T];X^{s-km}) ,\qquad
	\vert u\vert_{X^{s}_{(j)}}:=
	\sum_{k=0}^j \vert (\eps\dt)^ku\vert_{X^{s-km}_T},\\
	\label{Fs}
	F^s&:=&C([0,T];X^s)\times X^{s+m},\qquad \vert (f,g)\vert_{F^s}
	:=\vert f\vert_{X^s_T}+\vert g\vert_{s+m}
\end{eqnarray}
and, for all $(f,g)\in F^s$ and $t\in [0,T]$,
\begin{equation}
	\label{defI}
	{\mathcal I}^s(t,f,g):=\vert g\vert_s+\int_0^t \sup_{0\leq t''\leq t'}\vert f(t'')\vert_s dt'.
\end{equation}
\begin{assumption}
	\label{assintro3}
	Let $s_0,m $ and $T$ be as in Assumption \ref{assintro2}. There
	exist
	$d_1,d_1'\geq 0$
	such that for all $s\geq s_0+m$,
	$u^\eps\in X^{s+d_1}_{(1)}$
	and  $(f^\eps,g^\eps)\in F^{s+d_1'}$,
	the IVP
	\begin{equation}
	\label{assIVP}
	\dsp \dt \ve +\frac{1}{\eps}\cLe(t) \ve+\cFeu[t,\ue]\ve=f^\eps,
	\qquad
	\ve\init=g^\eps,
	\end{equation}
	admits a unique solution 
	$\ve\in C([0,T];X^{s})$ for all $\rge$, and 
	\begin{eqnarray*}
	\vert \ve\vert_{X^s_T}&\leq&
	C(\eps_0,s,T,\vert\ue\vert_{X^{s_0+m +d_1}_{(1)}})\\
	&\times&\big({\mathcal I}^{s+d_1'}(t,f^\eps,g^\eps)
	+\vert \ue\vert_{X^{s+d_1}_{(1)}}
	{\mathcal I}^{s_0+m+d_1'}(t,f^\eps,g^\eps)
	\big).
	\end{eqnarray*}
\end{assumption}
\begin{remark}
	The above energy estimate exhibits a loss of
	$d_1$ derivatives with respect to the reference state $\ue$
	(and of $d_1'$ derivatives with respect to the source term and
	initial data) in the sense that a control of $\ve$ in 
	$X^s_T$ requires a control of $\ue$ in $X^{s+d_1}_T$. This loss
	of information makes a standard Picard iterative scheme useless
	to find a solution to (\ref{intro1}). However, since the energy
	estimate is \emph{tame}, one can perform a Nash-Moser type
	iterative scheme. The fact that the energy estimate is also
	uniform with respect to $\rge$ is essential to obtain
	an existence time $\uT$ \emph{independent} of $\eps$.
\end{remark}

 \section{A Nash-Moser type theorem}\label{sectNM}

\subsection{Statement of the theorem}\label{sectstatement}

We state here the main theorem of this article (a generalization is also
given in Theorem \ref{th1b}' below). In the following statement, we use
the notations
$$
	\delta:=\max\{d_1,d_1'+m\},
	\qquad
	q:=D-m-d_1'
	\quad \mbox{ and }\quad
	P_{min}:=\delta+\frac{D}{q}\big(\sqrt{\delta}+\sqrt{2(\delta+q)}
	\big)^2,
$$
and we also recall that $F^{s+P}=C([0,T];X^{s+P})\times X^{s+P+m}$.
\begin{theorem}\label{th1}
	Let $T>0$, $s_0$, $m$, $d_1$ and $d_1'$ be such that
	Assumptions \ref{assintro1}-\ref{assintro3}
	are satisfied. Let also $D>\delta$,
	$P>P_{min}$,
	$s\geq s_0+m$ and
	$(h^\eps,\uinit)_\Rge$ be  bounded in
	$F^{s+P}$.\\
	Then there exists $0<\underline{T}\leq T$
	and a unique family 
	$(\underline{u}^\eps)_{\Rge}$  bounded in 
	$C([0,\underline{T}];X^{s+D})$ 
	and solving the IVPs (\ref{intro1})$_\Rge$.
\end{theorem}

\subsection{Proof of the theorem}\label{sectproof}

With the evolution operator $U^\eps(\cdot)$ defined in 
Assumption \ref{assintro1},
one can define
a nonlinear operator $\cGe[t,\cdot]$ as
$$
	\forall t\in [-T,T],\quad\forall u\in X^{s_0+m},\qquad
	\cGe[t,u]:=U^\eps(-t)\cFe[t,U^\eps(t)u].
$$
The next lemma shows that one can reduce the study of (\ref{intro1}) 
to the study of
\begin{equation}
	\label{eqalt}
	\left\lbrace
	\begin{array}{l}
	\dsp \dt \tue +\cGe[t,\tue]=U^\eps(-t)h^\eps\\
	\tue\init=\uinit
	\end{array}\right.
\end{equation}
and also states that $\cGe$ has the same properties as $\cFe$.
\begin{lemma}
	\label{lemproof1}
	{\bf i.} If $\tue\in C([0,T];X^s)$ solves (\ref{eqalt}) then 
	$\uue\in C([0,T];X^{s})$ solves (\ref{intro1}), where 
	$\uue(t):=U^\eps(t)\tue(t)$.\\
	{\bf ii.} Assumption \ref{assintro2} still holds if one replaces
	$\cFe$ by $\cGe$.
\end{lemma}
\begin{proof}
Assumption \ref{assintro1} shows that
if $\tue\in C([0,T];X^s)$ then 
$\uue\in C([0,T];X^{s})$. 
Remark now that if 
$\tue$ solves (\ref{eqalt}), then
$$
	\dt \big(U^\eps(t)\tue\big)=
	-\frac{1}{\eps}\cLe(t) U^\eps(t)\tue
	-U^\eps(t)\cGe[t,\tue]+h^\eps;
$$
since
$U^\eps(t)\cGe[t,\tue]=\cFe[t,U^\eps(t)\tue]$, 
the first point of the lemma follows. The second point 
is a direct consequence of 
Assumptions \ref{assintro1} and \ref{assintro2}.
\end{proof}

Defining  the space $F^s$ as in (\ref{Fs}) and 
$E^s$ as $C([0,T];X^s)\cap C^1([0,T];X^{s-m })$
endowed with its canonical norm (which makes $E^s$ different from $X^s_{(1)}$),
 we can use
Lemma \ref{lemproof1}, to check that finding a solution $\uue$ to
(\ref{intro1}) is  equivalent to finding a root $\tue$ of the equation
$\Phi^\eps(\tue)=0$, where
$$
	\Phi^\eps:
	\begin{array}{ccc}
	E^s & \to  & F^{s-m}\\
	u &\mapsto & (\underbrace{\dt u +\cGe[\cdot,u]-\wthe}_{:=\Phi_1(u)}, 
	u\init-\uinit),
	\end{array}
$$ 
for all $s\geq s_0+m$ and $\rge$, and with 
$\wthe(t):=U^\eps(-t)h^\eps(t)$.

We seek a root $\tue$ to the equation $\Phi^\eps(\tue)=0$ as the limit of 
a Nash-Moser type iterative
scheme, namely,
\begin{equation}
	\label{proof2}
	u_{k+1}^\eps=u_k^\eps+S_kv_k^\eps,
\end{equation}
with $S_k:={\mathcal S}_{\theta_k}$, for some $\theta_k>0$ to be determined,
and where $v_k^\eps$ solves
\begin{equation}
	\label{proof3}
	\left\lbrace
	\begin{array}{l}
	\dt v_k^\eps+\cGeu[t,u_k^\eps]v_k^\eps
	=-\Phi_1(u_k^\eps),\\
	v_k^\eps\,\init=\uinit-u_k^\eps\,\init.
	\end{array}\right.
\end{equation}
The following lemma shows that the above IVP can be solved and that the
knowledge of $u_k^\eps$ thus determines $v_k^\eps$.
\begin{lemma}
	\label{exist}
	Suppose that Assumptions \ref{assintro1}-\ref{assintro3} are satisfied,
	and let $s\geq s_0+m$. Assume also that  
	$u_k^\eps\in 
	E^{s+d_1}$ and 
	$\Phi^\eps(u_k^\eps)\in F^{s+d_1'}$.\\
	Then there exists a unique solution $v_k^\eps\in E^s$
	to (\ref{proof3}) and
	$$
	\vert v_k^\eps\vert_{E^s}\leq
	C(\eps_0,s,T,\vert u_k^\eps\vert_{E^{s_0+m+d_1}})
	\big(\vert \Phi^\eps(u_k^\eps)\vert_{F^{s+d_1'}}
	+\vert \Phi^\eps(u_k^\eps)\vert_{F^{s_0+m+d_1'}}
	\vert u_k^\eps\vert_{E^{s+d_1}}
	\big).
$$
\end{lemma}
\begin{proof}
From Assumption \ref{assintro3}, we know that there is a unique solution
$w_k^\eps$ of the IVP
$$
	\left\lbrace
	\begin{array}{l}
	\dt w_k^\eps+\frac{1}{\eps}\cLe(t)w_k^\eps+\cFeu[t,U^\eps(t)u_k^\eps]w_k^\eps=-U^\eps(t)\Phi_1(u_k^\eps),\\
	w_k^\eps\,\init=\uinit-u_k^\eps\,\init;
	\end{array}\right.
$$
as in the proof of Lemma \ref{lemproof1}, it is easy to check that 
$v_k^\eps:=U^\eps(-t)w_k^\eps$ solves (\ref{proof3}). \\
Since Assumption \ref{assintro1} implies that   
$\vert U^\eps(\cdot)u_k^\eps\vert_{X^{r}_{(1)}}
\leq C(\eps_0)\vert u_k^\eps\vert_{E^r}$ ($r\geq s_0+m$), one can deduce from
the estimate of Assumption \ref{assintro3} and
Assumption \ref{assintro1} that
\begin{equation}
	\label{controlr}
	\vert v_k^\eps\vert_{X^r_T}\leq
	C(\eps_0,r,T,\vert u_k^\eps\vert_{E^{s_0+m+d_1}})
	\big(\vert \Phi^\eps(u_k^\eps)\vert_{F^{r+d_1'}}
	+\vert \Phi^\eps(u_k^\eps)\vert_{F^{s_0+m+d_1'}}
	\vert u_k^\eps\vert_{E^{r+d_1}}
	\big);
\end{equation}
with $r=s$, this is the control we need on $\vert v_k^\eps\vert_{X^s_T}$;
to conclude the proof, we must therefore
show that the same bound holds for $\vert \dt v_k^\eps\vert_{X^{s-m}_T}$.
From the equation one has $\dt v_k^\eps=-\cGeu[t,u_k^\eps] v_k^\eps-\Phi_1(u_k^\eps)$, so that using Lemma \ref{lemproof1}.ii, one gets 
$$
	\vert \dt v_k^\eps\vert_{X^{s-m}_T}\leq
	C(s,T,\vert u_k^\eps\vert_{X^{s_0+m}_T} )
	\big(\vert v_k^\eps\vert_{X^s_T}
	+\vert v_k^\eps\vert_{X^{s_0+m}_T}\vert u_k^\eps\vert_{X^{s}_T}\big)
	+\vert \Phi_1(u_k^\eps)\vert_{X^{s-m}_T}.
$$
and one can conclude with (\ref{controlr}) (with $r=s$ and $r=s_0+m$).
\end{proof}
Let us now state the three lemmas which form the heart of the proof, 
and whose proof is postponed to the next subsections
for the sake of clarity.
\begin{lemma}
	\label{prooflem2}
	Let $D\geq m +d_1'$ and $s\geq s_0+m$.
	If, for some $M>0$, one has 
	$\vert u_j^\eps\vert_{E^{s+D}}\leq M$ ($j=k,k+1$),
	then
	$$
	\vert \Phi^\eps(u_{k+1}^\eps)\vert_{F^{s+d_1'}}\leq 
	C(s,T,M)\big( \theta_k^{m +d_1'-D}
	+\vert v_k^\eps\vert_{E^{s+D}}\big)
	\vert v_k^\eps\vert_{E^{s+D}},
	$$ 	
	with $C(s,T,M)$ independent of $\eps$.
\end{lemma}
\begin{lemma}
	\label{prooflem3}
	Let $D\geq d_1$ and $s\geq s_0+m$.
	If, for some $M>0$, one has 
	$\vert u_{k+1}^\eps\vert_{E^{s+D}}\leq M$,
	then
	$$
	\vert v_{k+1}^\eps\vert_{E^{s}}\leq
	C(\eps_0,s,T,M) \vert \Phi^\eps(u_{k+1}^\eps)\vert_{F^{s+d_1'}},
	$$ 	
	with $C(\eps_0,s,T,M)$ independent of $\eps$.
\end{lemma}
\begin{lemma}
	\label{prooflem4}
	Let $\delta:=\max\{d_1,(d_1'+m)\}$, $P\geq D\geq \delta$ and  
	$s\geq s_0+m$.
	If, for some $M>0$, one has 
	$\vert u_{k}^\eps\vert_{E^{s+D}}\leq M$
	and $\vert (h^\eps,\uinit)\vert_{F^{s+P-m}}\leq M$, then
	$$
	\vert u_{k+1}^\eps\vert_{E^{s+P}}\leq
	C(\eps_0,s,T,M)(1+\theta_k^\delta)
	(1+\vert u_k^\eps\vert_{E^{s+P}}).
	$$
	If moreover $\vert u_{k+1}^\eps\vert_{E^{s+D}}\leq M$, then one also
	has
	$$
	\vert v_{k+1}^\eps\vert_{E^{s+P-\delta}}
	\leq C(\eps_0,s,T,M)(1
	+\vert u_{k+1}^\eps\vert_{E^{s+P}}).
	$$
\end{lemma}

\bigbreak

We can 
now proceed with the proof of the theorem, which is a typical 
Nash-Moser iterative scheme : Lemmas \ref{prooflem2} 
and \ref{prooflem3} provide a control of $\vert v_{k+1}^\eps\vert_{E^s}$ 
in terms of 
$\vert v_k^\eps\vert_{E^{s+D}}$, thus exhibiting a loss of $D$ derivatives 
but providing  a 
rapid decay of $\vert v_{k+1}^\eps\vert_{E^s}$,
while Lemma \ref{prooflem4} control the growth of 
$\vert v_{k+1}^\eps\vert_{E^{s+P-\delta}}$.
A control of $\vert v_{k+1}^\eps\vert_{E^{s+D}}$ is then recovered
by the interpolation
formula
\begin{equation}
	\label{interpform}
	\vert v_{k+1}^\eps\vert_{E^{s+D}}\leq\cst
	\vert v_{k+1}^\eps\vert_{E^s}^\mu\,
	 \vert v_{k+1}^\eps\vert_{E^{s+P-\delta}}^{1-\mu},
\end{equation}
with $\mu=1-\frac{D}{P-\delta}$.

Before entering the heart of the proof, let us define the
sequence $(\theta_k)_{k}$ used for the smoothing operators as
$\theta_{k+1}=\theta_k^r$ ($k\in\N$), for some $r>1$ defined
below.
\begin{remark}
	\label{remproof1}
	One has $\theta_k=\theta_0^{r^k}$, so that
	(if $r>1$),
	$\sum_{k\in\N}\theta_k^{-q}=:\underline{\theta}$ converges if 
	and only if $\theta_0>1$. 
	Moreover, $\underline{\theta}$ can be 
	made arbitrarily small provided that $\theta_0$ is chosen 
	large enough. 
\end{remark}

We are now set to control the sequences $(u_k^\eps)_{k\in \N}$ and
$(v_k^\eps)_{k\in \N}$ by induction. For some $M>0$  such that
\begin{equation}
	\label{condM}
	\vert (h^\eps,\uinit)\vert_{F^{s+P-m}}\leq M
	\quad\mbox{ and }\quad
	\vert u_0^\eps\vert_{E^{s+D}}\leq M/2,
\end{equation}
we define the properties (i)\sk-(iii)\sk  as
\begin{itemize}
	\item (i)\sk: $\vert u_{k}^\eps\vert_{E^{s+P}}\leq  
	\theta_k^\alpha$;
	\item (ii)\sk: $\vert u^\eps_{k}\vert_{E^{s+D}}\leq M$;
	\item (iii)\sk: $\vert v_k^\eps\vert_{E^{s+D}}
	\leq \theta_k^{-q}$, with $q=D-m -d_1'>0$.
\end{itemize}

\noindent
\emph{Proof of (i)$_{k+1}$-(iii)$_{k+1}$ assuming (i)\sk-(iii)\sk.} Since
one has $\vert u_k^\eps\vert_{E^{s+D}}\leq M$ by (ii)\sk, 
$\vert u_{k}^\eps\vert_{E^{s+P}}\leq \theta_k^\alpha$ by (i)\sk  and 
$\vert (h^\eps,\uinit)\vert_{F^{s+P-m}}\leq M$ by definition of $M$, 
one can apply Lemma \ref{prooflem4} to obtain
\begin{eqnarray*}
	\vert u_{k+1}^\eps\vert_{E^{s+P}}
	&\leq& C(\eps_0,s,T,M)(1+\theta_k^\delta)(1+\theta_k^\alpha)\\
	&=& f(\eps_0,s,T,M,k)\theta_{k+1}^\alpha,
\end{eqnarray*}
with
$$
	f(\eps_0,s,T,M,k)=C(\eps_0,s,T,M)(1+\theta_k^\delta)
	(1+\theta_k^\alpha)
	\theta_k^{-\alpha r}.
$$
Assuming that
\begin{equation}
	\label{cond}  
	\delta-\alpha(r-1)<0,
\end{equation}
it follows from the
explicit expression of $f(\eps_0,s,T,M,k)$ that $f(\eps_0,s,T,M,k)\leq 1$ 
for all $k\in\N$ provided that $\theta_0$ is chosen
large enough. This proves (i)\skp.\\
Recalling that $u_{k+1}^\eps=u_k^\eps+S_kv_k^\eps$, one has 
$u_{k+1}^\eps=u_0^\eps+\sum_{j=0}^k S_jv_j^\eps$, and thus
$\vert u_{k+1}^\eps\vert_{E^{s+D}}\leq M/2+\sum_{j=0}^k \theta_k^{-q}$. 
As seen in Remark \ref{remproof1}, one then gets  (ii)\skp provided that
$\theta_0$ is chosen large enough.\\
In order to prove (iii)\skp, remark first that it follows from
Lemmas \ref{prooflem2} and \ref{prooflem3} and the choice of the sequence
$(\theta_k)_{k\in\N}$ that
\begin{equation}
	\label{proof5}
	\vert v_{k+1}^\eps\vert_{E^s}\leq C(\eps_0,s,T,M)\theta_{k+1}^{-2q/r}.
\end{equation}
We can also use the second
assertion of Lemma \ref{prooflem4} to obtain
\begin{equation}
	\label{proof6}
	\vert v_{k+1}^\eps\vert_{E^{s+P-\delta}}
	\leq C(\eps_0,s,T,M)(1+\theta_{k+1}^\alpha).
\end{equation}
It follows therefore from (\ref{interpform}), (\ref{proof5}) 
and (\ref{proof6}) that
\begin{eqnarray*}
	\vert v_{k+1}^\eps\vert_{E^{s+D}}
	&\leq& 
	C(\eps_0,s,T,M) \theta_{k+1}^{-2\mu q/r}(1+\theta_{k+1}^\alpha)^{1-\mu}
\\
	&=& g(\eps_0,s,T,M,k)\theta_{k+1}^{-q},
\end{eqnarray*}
with $\dsp g(\eps_0,s,T,M,k):=C(\eps_0,s,T,M) 
\theta_{k+1}^{-2\mu q/r}(1+\theta_{k+1}^\alpha)^{1-\mu}\theta_{k+1}^q$. 
Choosing $r$ such that
\begin{equation}
	\label{choice}
	1<r<\frac{2\mu q}{q+\alpha(1-\mu)},
\end{equation} one
gets that $g(\eps_0,s,T,M,k)\leq 1$ for all $k\in \N$, provided that $\theta_0$
is chosen large enough.\\
It follows from the lines above that in 
order to complete the proof of the heredity of the induction property,
we just have to take $\theta_0$ large enough, and to prove that
one can choose $\alpha$, $r$ and $P$ 
such that the conditions
(\ref{cond}) and  (\ref{choice}) are satisfied. 
This is done in the following lemma:
\begin{lemma}
	Let $\alpha= \delta+\sqrt{2\delta(\delta+q)}$; if  
	$P>\delta+\frac{D}{q}(\sqrt{\delta}+\sqrt{2(\delta+q)})^2$, 
	there exists $r>1$ such that conditions (\ref{cond}) 
	and 
	(\ref{choice}) are satisfied.
\end{lemma}
\begin{proof}
Let us denote ${\underline r}:=\frac{2\mu q}{q+\alpha(1-\mu)}$. Quite
obviously, (\ref{cond}) and (\ref{choice}) are satisfied with 
$r=\underline{r}-\epsilon$ ($\epsilon>0$ small enough), provided that 
$\underline{r}-1>\delta/\alpha$,
that is,
$$
	\frac{(2q+\alpha)\mu-(q+\alpha)}{q+\alpha(1-\mu)}>\frac{\delta}{\alpha},
$$
or equivalently, if
$$
	\mu>1-\frac{q(1-\delta/\alpha)}{2q+\alpha+\delta}=:\mu_{\min}(\alpha).
$$
The value of $\alpha$ given in the statement of the lemma corresponds to
the minimum of $\mu_{min}(\alpha)$. One then computes that 
$\mu_{min}(\alpha)=1-\frac{q}{(\sqrt{\delta}+\sqrt{2(\delta+q)})^2}$, and the lemma then follows from the observation
that $\mu>\mu_{min}(\alpha)$ is equivalent to $P>\delta+\frac{D}{1-\mu_{min}}$.
\end{proof}

\noindent
\emph{Proof of (i)$_0$-(iii)$_0$.} We have to construct here the first 
term of the sequence $u_0^\eps$ in such a way that (i)$_0$-(iii)$_0$ and
(\ref{condM})
are satisfied for some $M>0$ and $\theta_0>0$. We need the following lemma:
\begin{lemma}\label{lmmil}
	For all $s\geq s_0+m$ and $(h^\eps,\uinit)\in F^{s+P}$, 
	there exists $u_0^\eps\in E^{s+P}$
	such that $u_0^\eps\,\init=\uinit$ and such that
	$$
	\vert u_0^\eps\vert_{E^{s+D}}\leq 
	C(s,T,\vert (h^\eps,\uinit)\vert_{F^{s+D}})
	\quad\mbox{ and }\quad
	\vert u_0^\eps\vert_{E^{s+P}}\leq 
	C(s,T,\vert (h^\eps,\uinit)\vert_{F^{s+P}}),
	$$
	and
	$$
	\vert\Phi^\eps(u_0^\eps)\vert_{F^{s+D+d_1'}}	
	\leq TC(s,T,\vert (h^\eps,\uinit)\vert_{F^{s+D+d_1'+m}}).
	$$
\end{lemma}
\begin{proof}
Let us define $u_0^\eps\in C([0,T];X^{s+P})$ as
$$
	u_0^\eps(t)=\uinit+\int_0^t\big(\wthe(t')-\cGe[t',\uinit]\big)dt'.
$$
From Lemma \ref{lemproof1} and the definition of 
$\vert\cdot\vert_{E^s}$, one gets
for all $r\geq 0$,
\begin{equation}
	\label{unedeplus}
	\vert u_0^\eps\vert_{E^{s+r}}\leq
	\vert \uinit\vert_{s+r}+C(s,T,\vert \uinit \vert_{s_0+m})
	\big(\vert h^\eps\vert_{X^{s+r}_T}
	+\vert \uinit\vert_{s+r+m}\big);
\end{equation}
the estimates on $u_0^\eps$ given in the lemma are thus a consequence of
(\ref{unedeplus}), with $r=D$ and $r=P$.\\
By definition of $\Phi^\eps$, one also has
\begin{eqnarray*}
	\Phi^\eps(u_0^\eps)&=&
	\big(\cGe[\cdot,u_0^\eps]-\cGe[\cdot,\uinit],0\big)\\
	&=&\big(\int_0^1 \cGeu[\cdot,\uinit+z(u_0^\eps-\uinit)](u_0^\eps-\uinit)dz,0\big),
\end{eqnarray*}
so that one deduces from Assumptions \ref{assintro1} and \ref{assintro2} that
$$
	\vert \Phi^\eps(u_0^\eps)\vert_{F^{s+D+d_1'}}\leq
	C(s,T,\vert u_0^\eps\vert_{X^{s+D+d_1'+m}_T},
	\vert \uinit\vert_{s+D+d_1'+m})
	\vert u_0^\eps-\uinit\vert_{X^{s+D+d_1'+m}_T},
$$
and the estimate on $\Phi^\eps(u_0^\eps)$ of the lemma follows easily.
\end{proof}
Thanks to the lemma, taking 
$M=M(s,T,\vert(h^\eps,\uinit)\vert_{F^{s+D}})$ large enough, one
gets $\vert u_0^\eps\vert_{E^{s+D}}\leq M/2$, which proves (ii)$_0$. Choosing
$\theta_0=\theta_0(s,T,\vert (h^\eps,\uinit)\vert_{F^{s+P}})$ large enough, one
also gets (i)$_0$ from Lemma \ref{lmmil}. 
In order to prove (iii)$_0$, remark first
that Lemma \ref{exist} yields $\vert v_0^\eps\vert_{E^{s+D}}\leq 
C(\eps_0,s,T,M)\vert 
\Phi^\eps(u_0^\eps)\vert_{F^{s+D+d_1'}}(1+\theta_0^\alpha)$. 
It follows therefore from the lemma that,
taking a smaller $T$ if necessary, (iii)$_0$ is satisfied, which ends
the induction proof of properties (i)$_k$, (ii)$_k$ and (iii)$_k$.

\bigbreak

The end of the existence part of the proof of the theorem is now 
straightforward: it follows from
(i)$_k$, (ii)$_k$ and (iii)$_k$ that the series $u_0^\eps+\sum_k S_k v_k^\eps$
converges to some $u^\eps\in E^{s+D}$ and taking the limit $k\to\infty$
in Lemma \ref{prooflem2} shows that $\Phi^\eps(u^\eps)=0$.

In order to conclude the proof of the theorem, we must now prove that
the solution constructed above is unique. Assuming that $u^{\eps,j}\in E^{s+D}$
($j=1,2$) are both solutions to (\ref{intro1}), 
we show that $w:=u^{\eps,2}-u^{\eps,1}$ is
identically $0$. Let us remark that $w$ solves the IVP
$$
	\left\lbrace
	\begin{array}{l}
	\dt w+\frac{1}{\eps}\cLe(t)w+\cFeu[t,u^{\eps,2}]w=H,\\
	w\init=0,
	\end{array}\right.
$$
with $H:=\cFe[t,u^{\eps,1}]-\cFe[t,u^{\eps,2}]-\cFeu[t,u^{\eps,2}](u^{\eps,1}-u^{\eps,2})$.\\
A direct application of Assumption \ref{assintro3} yields
$$
	\vert w(t)\vert_{s_0+m}\leq C(\eps_0,T,\vert u^{\eps,1}\vert_{X^{s_0+m+\delta}_{(1)}})\int_0^t
	\sup_{0\leq t''\leq t'}\vert H(t'')\vert_{s_0+m+d_1'} dt',
$$
and since
$\vert H(t')\vert_{s_0+m+d_1'}\leq C(s,T, 
\vert u^{\eps,2}\vert_{X^{s_0+m+\delta}_T})\vert w(t')\vert_{s_0+m}$ by
Assumption \ref{assintro2}$_{(\ref{condass3})}$, a Gronwall argument 
shows that $w$=0.

\subsection{Proof of Lemmas \ref{prooflem2}, \ref{prooflem3}, \ref{prooflem4}}
\label{sectprooflemmas}

\subsubsection{Proof of Lemma \ref{prooflem2}}

In order to give an upper bound for $\vert \Phi^\eps(u_{k+1}^\eps)\vert_{F^{s+d_1'}}$, we need to control $\Phi_1(u_{k+1}^\eps)$ in $X^{s+d_1'}_T$ and
$\vert u_{k+1}^\eps \,\init-\uinit\vert_{s+d_1'+m}$.\\
First remark that a second order Taylor expansion of $\Phi_1(u_{k+1}^\eps)$
yields
\begin{eqnarray*}
	\Phi_1(u_{k+1}^\eps)&=&\Phi_1(u_k^\eps)+\Phi_1'(u_k^\eps)(u_{k+1}^\eps-u_k^\eps)\\
	& &+\int_0^1(1-z)\Phi_1''(u_k^\eps+z(u_{k+1}^\eps-u_k^\eps))
	(u_{k+1}^\eps-u_k^\eps,u_{k+1}^\eps-u_k^\eps)dz.
\end{eqnarray*}
Since by (\ref{proof2}), one has $u_{k+1}^\eps-\ue_k=S_k v_k^\eps$, and
since by definition
$\Phi_1(u_{k}^\eps)=\dt u_{k}^\eps+\cGe[\cdot,u_{k}^\eps]-\wthe$, it
follows that
\begin{equation}
	\label{proofl1}
	\Phi_1(u_{k+1}^\eps)=E_1+E_2,
\end{equation}
with 
\begin{eqnarray}
	\nonumber
	E_1&=&\Phi_1(u_k^\eps)+\dt v_k^\eps
	+\cGeu[\cdot,u_k^\eps]v_k^\eps\\
	\nonumber
	& &+\int_0^1(1-z)\cGeuu[\cdot,u_k^\eps+z (u_{k+1}^\eps-u_k^\eps)](S_kv_k^\eps,S_k v_k^\eps )dz\\
	\label{eqE1}
	&=& \int_0^1(1-z)\cGeuu[\cdot,u_k^\eps+z (u_{k+1}^\eps-u_k^\eps)](S_kv_k^\eps,S_k v_k^\eps )dz
\end{eqnarray}
(the last equality stemming from the fact that $v_k^\eps$ solves (\ref{proof3})),
and
\begin{equation}
	\label{eqE2}
	E_2=(S_k-1)\dt v_k^\eps
	+\cGeu[\cdot,u_k^\eps]\big((S_k-1)v_k^\eps\big).
\end{equation}
Since by Lemma \ref{lemproof1}.ii, $\cGe$ satisfies Assumption \ref{assintro2}$_{(\ref{condass3})}$, 
and since $s+d_1'+m\leq s+D$,
one can control $E_1$ as
\begin{eqnarray}
	\nonumber
	\vert E_1\vert_{X^{s+d_1'}_T}
	&\leq& C(s,T,\vert u_k^\eps\vert_{X^{s+D}_T}, 
	\vert u_{k+1}^\eps\vert_{X^{s+D}_T})
	\vert S_k v_k^\eps\vert^2_{X^{s+D}_T}\\
	\label{proofl2}
	&\leq&C(s,T,M) \vert v_k^\eps\vert_{X^{s+D}_T}^2.
\end{eqnarray}
Since Lemma \ref{lemproof1}.ii also ensures that $\cGe$ satisfies
Assumption \ref{assintro2}$_{(\ref{condass2})}$, one gets 
$$
	\vert \cGeu[\cdot,u_k^\eps]\big((S_k-1)v_k^\eps\big)
	\vert_{X^{s+d_1'}_T}
	\leq C(s,T,M)
	\supT\vert (S_k-1)v_k^\eps(t)\vert_{s+m +d_1'}.
$$
It is then a consequence of the properties of the regularizing operators 
(recall that $S_k={\mathcal S}_{\theta_k}$), that
\begin{eqnarray}
	\nonumber
	\vert E_2\vert_{X^{s+d_1'}_T}&\leq& \cst \theta_k^{m +d_1'-D}\big(
	\vert\dt v_k^\eps\vert_{X^{s+D-m}_T }
	+C(s,T,M)\vert v_k^\eps\vert_{X^{s+D}_T}\big)\\
	\label{proofl3}
	&\leq& C(s,T,M)\theta_k^{m +d_1'-D}\vert v_k^\eps\vert_{E^{s+D}}.
\end{eqnarray}
It is then a simple consequence of (\ref{proofl1}), (\ref{proofl2})
and (\ref{proofl3}) to conclude that
\begin{equation}
	\label{proofl4}
	\vert \Phi_1(u_{k+1}^\eps)\vert_{X^{s+d_1'}_T}\leq 
	C(s,T,M)\big( \theta_k^{m+d_1'-D}+\vert v_k^\eps\vert_{E^{s+D}}\big)
	\vert v_k^\eps\vert_{E^{s+D}}.
\end{equation}
We now turn to control
$\vert u_{k+1}^\eps \,\init-\uinit\vert_{s+d_1'+m}$. Since 
$u_{k+1}^\eps\,\init-\uinit=(S_k-1)v_k^\eps\,\init$,
one gets
\begin{eqnarray}
	\nonumber
	\vert u_{k+1}^\eps \,\init-\uinit\vert_{s+d_1'+m}&\leq& 
	\cst \theta_k^{m+d_1'-D}\supT \vert v_k^\eps(t)\vert_{s+D}\\
	\label{proofl5}
	&\leq &
	\cst \theta_k^{m+d_1'-D}\vert v_k^\eps\vert_{E^{s+D}}.
\end{eqnarray}
The lemma follows directly from (\ref{proofl4}) and (\ref{proofl5}).

\subsubsection{Proof of Lemma \ref{prooflem3}}

Since
$\vert u_{k+1}^\eps\vert_{E^{s +d_1}}\leq M$,
 one gets therefore from Lemma \ref{exist} (at step $k+1$),
\begin{equation}
	\label{proofl6}
	\vert v_{k+1}^\eps\vert_{E^s}\leq
	C(\eps_0,s,T,M)
	\vert \Phi^\eps(u_{k+1}^\eps)\vert_{F^{s+d_1'}},
\end{equation}
and the lemma is proved.

\subsubsection{Proof of Lemma \ref{prooflem4}}

Thanks to Lemma \ref{lemproof1}.ii, one has, 
for all $r\geq s_0$,
\begin{equation}
	\label{proofl8}
	\vert \Phi^\eps(u)\vert_{F^r}\leq
	C(\vert u\vert_{E^{s_0+m }})\vert u\vert_{E^{r+m }}
	+\cst \vert (h^\eps,\uinit)\vert_{F^r};
\end{equation}
remark also that since $u_{k+1}^\eps=u_k^\eps+S_k v_k^\eps$, one can use the properties
of the regularizing operator $S_k={\mathcal S}_{\theta_k}$ to obtain
\begin{equation}
	\label{proofl9}
	\vert u_{k+1}^\eps\vert_{E^{s+P}}\leq
	\vert u_{k}^\eps\vert_{E^{s+P}}
	+\cst \theta_k^{\delta}
	\vert v_k^\eps\vert_{E^{s+P-\delta}}.
\end{equation}
From Lemma \ref{exist},  one deduces
$$
	\vert v_k^\eps\vert_{E^{s+P-\delta}}\leq C(\eps_0,s,T,M)
	\Big( \vert \Phi^\eps(u_k^\eps)\vert_{F^{s_0+m +d_1'}}
	 \vert u_k^\eps\vert_{E^{s+P}}
	+\vert \Phi^\eps(u_k^\eps) \vert_{F^{s+P-m}}\Big)
$$
so that, using (\ref{proofl8}) with $r=s_0+m+d_1'$ and $r=s+P-m$, and
the assumption made on $(h^\eps,\uinit)$, one 
obtains
\begin{equation}
	\label{proofl10}
	\vert v_k^\eps\vert_{E^{s+P-\delta}}
	\leq
	C(\eps_0,s,T,M)(1
	+\vert u_k^\eps \vert_{E^{s+P}}).
\end{equation}
Together with (\ref{proofl9}), this last estimate shows that
$$
	\vert u_{k+1}^\eps\vert_{E^{s+P}}\leq
	C(\eps_0,s,T,M)(1+\theta_k^\delta)(1
	+\vert u_k^\eps\vert_{E^{s+P}}),
$$
so that the proof of the first assertion is complete.\\
The last part of the lemma is exactly (\ref{proofl10}) with the
index $k$ replaced by $k+1$.

\section{Further results}\label{sectfurther}

We propose in this section a more general version of Theorem \ref{th1} and
some remarks extending its range of validity. We also
a stability property very useful for 
the justification of asymptotic models
for instance.

\subsection{A more general version of Theorem \ref{th1}}\label{sectgen}

The aim of this section is to prove a result similar to Theorem \ref{th1}
when the energy estimates of Assumption \ref{assintro3} involve
$p+1$ ($p\geq 1$) time derivatives of the reference solution $\ue$
(such a situation occurs for instance with the water-waves equations). 
With this goal in mind,
we replace the three assumptions \ref{assintro1}-\ref{assintro3} by
generalizations to the case $p\neq 0$. We first generalize the
spaces $E^s$ and $F^s$ used in the proof of Theorem \ref{th1} as follows:
$$
	E^{s}_{(p+1)} :=\bigcap_{i=0}^{p+1} C^i([0,T];X^{s-im}),
	\qquad
	F^{s}_{(p)} := \bigcap_{i=0}^p C^i([0,T];X^{s-im})\times X^{s+m}
$$
endowed with the norms
\begin{eqnarray*}
	\vert u\vert_{E^{s}_{(p+1)}}=\vert u\vert_{X^s_T}
	+\vert \dt u\vert_{X^{s-m}_T}
	+\sum_{i=1}^p \vert (\eps\dt)^i \dt u\vert_{X^{s-(i+1)m}_T},\\
	\vert (f,g)\vert_{F^{s}_{(p)}}=\vert f\vert_{X^s_T}+\vert g\vert_{s+m}
	+\sum_{i=1}^p \vert (\eps\dt)^{i} f\vert_{X^{s-im}_T}
\end{eqnarray*}
and we also define for all $(f,g)\in F^s_{(p)}$ and $t\in [0,T]$,
$$
	{\mathcal I}^{s}_{(p)}(t,f,g)=\vert g\vert_{s+m}+\int_0^t
	\sum_{i=0}^p \sup_{0\leq t''\leq t'}
	\vert (\eps\dt)^i f(t'')\vert_{s-im}dt'
$$
(so that $E^{s}_{(1)}$, $F^{s}_{(0)}$ and ${\mathcal I}^s_{(0)}$ coincide with 
$E^s$, $F^s$ and ${\mathcal I}^s$ respectively).
\begin{assumptionbis}\label{assb1}
	Assumption \ref{assintro1}
	holds with (\ref{condass11}) replaced by the stronger condition
	(when $p\geq 1$):
	\begin{enumerate}
		\item[(\ref{condass11})'] For all $s\geq s_0$, one has
	$\cLe\in C^p(\R;{\mathfrak L}(X^{s+m};X^{s}))$ and 
	for all $i=0,\dots,p$, 
	$(\eps^i\frac{d^i}{(dt)^i}\cLe(\cdot))_\Rge$ is 
	bounded in $C([0,T];{\mathfrak L}(X^{s+m};X^{s}))$.
	\end{enumerate}
\end{assumptionbis}
\begin{assumptionbis}\label{assb2}
	Assumption \ref{assintro2} 
	holds with (\ref{condass1})-(\ref{condass3}) replaced replaced
	by the stronger conditions:
	For all 
	$0\leq i\leq p$ and $0\leq i+j\leq p+2$, and for all $s\geq s_0+im$, 
	one has $\cFe\in C^i([0,T];C^j(X^{s+m},X^{s-im}))$ and
	\begin{eqnarray*}
	\lefteqn{\supT\vert \eps^i {\mathcal F}^{\eps(i)}_{(j)}[t,u]
	(v_1,\dots,v_{j})
	\vert_{s-im}\leq
	C(s,T,\vert u\vert_{s_0+(i+1)m })}\\
	& &\times\big(\sum_{k=1}^{j} \vert v_k\vert_{s+m}
	\prod_{l\neq k }\vert v_l\vert_{s_0+(i+1)m}
	+\vert u\vert_{s+m}
	\prod_{k=1}^{j} \vert v_k\vert_{s_0+(i+1)m}\big).
	\end{eqnarray*}
\end{assumptionbis}
\begin{assumptionbis}
	\label{assb3}
	There exists $p\in\N$ such that for all $s\geq s_0+m$,
	$u^\eps\in X^{s+d_1}_{(p+1)}$
	and  $(f^\eps,g^\eps)\in F^{s+d_1'}_{(p)}$,
	the IVPs (\ref{assIVP})$_\Rge$ admit a unique solution
	$\ve\in C([0,T];X^{s})$, and
	\begin{eqnarray*}
	\forall t\in [0,T],\qquad
	\vert \ve\vert_{X^s_T}&\leq&
	C(\eps_0,s,T,\vert\ue\vert_{X^{s_0+m +d_1}_{(p+1)}})\\
	&\times &\big({\mathcal I}_{(p)}^{s+d_1'}(t,f^\eps,g^\eps)
	+\vert \ue\vert_{X^{s+d_1}_{(p+1)}}
	{\mathcal I}_{(p)}^{s_0+m+d_1'}(t,f^\eps,g^\eps)
	\big).
	\end{eqnarray*}
\end{assumptionbis}

Theorem \ref{th1} then admits the following generalization
(with $\delta$ and $P_{min}$ as defined in Theorem \ref{th1}):
\begin{theorembis}\label{th1b}
	Let $p\in\N$, $T>0$, $s_0$, $m$, $d_1$ and $d_1'$ be such that
	Assumptions \ref{assb1}'-\ref{assb3}'
	are satisfied. Let also $D>\delta$,
	$P>P_{min}$,
	 $s\geq s_0+(p+1)m$
	and $(h^\eps,\uinit)_\Rge$ be  bounded in
	$F^{s+P}_{(p)}$.\\
	Then there exists $0<\underline{T}\leq T$
	and a unique family 
	$(\underline{u}^\eps)_{\Rge}$  bounded in 
	$C([0,\underline{T}];X^{s+D})$ 
	and solving the IVPs (\ref{intro1})$_\Rge$.
\end{theorembis}
\begin{proof}
The proof is a generalization of the proof of Theorem \ref{th1},
and we just sketch the adaptations to be done.\\

The second property of Lemma \ref{lemproof1} can be generalized as
follows:
\begin{lemmabis}\label{lmb1}
	Assumption \ref{assb2}' still holds if one replaces $\cFe$ by
	$\cGe$.
\end{lemmabis}
\begin{proof}
Let us first prove the following fact (recalling
that the normed space $X^{s}_{(i)}$ is defined in (\ref{Xs})):
for $0\leq i\leq p$ and $s\geq s_0+im$,
\begin{equation}
	\label{fact}
	\forall f\in X^{s}_{(i)},\quad U^\eps(\cdot)f\in X^{s}_{(i)}
	\quad\mbox{ and }\quad
	\vert U^\eps(\cdot)f\vert_{X^{s}_{(i)}}
	\leq C(s,T)\vert f\vert_{X^{s}_{(i)}};
\end{equation}
indeed, from the definition of the
evolution operator $U^\eps(\cdot)$, one can check that 
$(\eps\dt)^i \big(U^\eps(\cdot)f(\cdot)\big)$ is a sum of terms
of the form
$$
	(\cLe)^{\alpha_0}(\eps\dt \cLe)^{\alpha_1}\dots
	(\eps^{i-1}\dt^{i-1}\cLe)^{\alpha_{i-1}}U^\eps(\cdot)
	\big((\eps\dt)^\beta
	f\big),
$$
with $\alpha_0+2\alpha_1+\dots+i\alpha_{i-1}+\beta=i$, so that
(\ref{fact}) is a direct consequence of Assumption
\ref{assb1}'.\\
From the definition of $\cGe$, one has, for all $0\leq i+j\leq p+2$,
$$
	\eps^i {\mathcal G}^{\eps(i)}_{(j)}[\cdot,u](v_1,\dots,v_{j})
	=(\eps\dt)^i \big(U^\eps(-\cdot)
	f(\cdot)\big),
$$
with $f(t)={\mathcal F}_{(j)}^\eps[t,U^\eps(t)u]
(U^\eps(t)v_1,\dots,U^\eps(t)v_{j})$ so that one 
deduces directly from (\ref{fact}) that for all $s\geq s_0+im$,
\begin{equation}
	\label{step1}
	\vert (\eps\dt)^i {\mathcal G}^\eps_{(j)}[\cdot,u](v_1,\dots,v_{j})
	\vert_{X^{s-im}_T}\leq C(s,T)\vert f\vert_{X^{s}_{(i)}}.
\end{equation}
Writing $\uu:=U^\eps(\cdot)u$ and $\uv_q:=U^\eps(\cdot)v_q$, one has
$f={\mathcal F}^\eps_{(j)}[t,\uu](\uv_1,\dots,\uv_{j})$ and for all 
$0\leq l\leq i$, $(\eps\dt)^l f$ is a sum of terms of the form
$$
	\eps^{l_0}{\mathcal F}^{\eps(l_0)}_{(j+\gamma_1+\dots+\gamma_l)}
	[t,\uu]((\eps\dt)^{l_1}\uv_1,\dots,(\eps\dt)^{l_{j}}\uv_{j},
	[\eps\dt \uu]^{\gamma_1},\dots,[(\eps\dt)^l \uu]^{\gamma_l}),
$$
with $l_0+l_1+\dots+l_{j}+\gamma_1+2\gamma_2+\dots+l\gamma_l=l$, and
where $[(\eps\dt)^l\uu]^{\gamma_l}$ stands for the $\gamma_l$-uplet with
$(\eps\dt)^l\uu$ on each component. It follows therefore from
Assumption \ref{assb2}' that
\begin{eqnarray}
	\nonumber
	\lefteqn{\vert (\eps\dt)^l f\vert_{X^{s-lm}_T}
	\leq C(s,T,\vert \uu\vert_{X^{s_0+(l+1)m}_{(l)}})}\\
	\label{step2}
	& &\times
	\big(
	\sum_{q=1}^{j}
	\vert \uv_q\vert_{X^{s+m}_{(l)}}\prod_{q'\neq q}
	\vert \uv_{q'}\vert_{X^{s_0+(l+1)m}_{(l)}}
	+\vert \uu\vert_{X^{s+m}_{(l)}}\prod_{q=1}^{j}
	\vert \uv_q\vert_{X^{s_0+(l+1)m}_{(l)}}\big)
\end{eqnarray}
and since (\ref{fact}) implies that for all $r\geq s_0+lm$,
$\vert \uu\vert_{X^{r}_{(l)}}\leq C(r,T)\vert u\vert_{X^r_T}=C(r,T)\vert u\vert_r$ and 
$\vert \uv_q\vert_{X^{r}_{(l)}}\leq C(r,T)\vert v_q \vert_{X^r_T}=
C(r,T)\vert v_q\vert_r$ 
($q=1,\dots,j$),
one deduces
\begin{eqnarray*}
	\lefteqn{\forall 0\leq l\leq i,\qquad
	\vert (\eps\dt)^l f\vert_{X^{s-lm}_T}
	\leq C(s,T,\vert u\vert_{{s_0+(l+1)m})}}\\
	& &\times
	\big(
	\sum_{q=1}^{j}
	\vert v_q\vert_{s+m}\prod_{q'\neq q}
	\vert v_{q'}\vert_{{s_0+(l+1)m}}
	+\vert u\vert_{{s+m}}\prod_{q=1}^{j}
	\vert v_q\vert_{{s_0+(l+1)m}}\big).
\end{eqnarray*}
Together with (\ref{step1}), this shows that $\cGe$ satisfies 
Assumption \ref{assb2}'.
\end{proof}
Lemma \ref{exist} can then be generalized as follows:
\begin{lemmabis}
	\label{existbis}
	Suppose that Assumptions \ref{assintro1}'-\ref{assintro3}'
	 are satisfied,
	and let $s\geq s_0+(p+1)m$. Assume also that  
	$u_k^\eps\in 
	E^{s+d_1}_{(p+1)}$ and 
	$\Phi^\eps(u_k^\eps)\in F^{s+d_1'}_{(p)}$.\\
	Then there exists a unique solution $v_k^\eps\in E^{s}_{(p+1)}$
	to (\ref{proof3}) and
	$$
	\vert v_k^\eps\vert_{E^{s}_{(p+1)}}\leq
	C(\eps_0,s,T,\vert u_k^\eps\vert_{E^{s_0+(p+1)m+d_1}_{(p+1)}})
	\big(\vert \Phi^\eps(u_k^\eps)\vert_{F^{s+d_1'}_{(p)}}
	+\vert \Phi^\eps(u_k^\eps)\vert_{F^{s_0+m+d_1'}_{(p)}}
	\vert u_k^\eps\vert_{E^{s+d_1}_{(p+1)}}
	\big).
	$$
\end{lemmabis}
\begin{proof}
Following the same steps as in the proof of Lemma \ref{exist}, one 
can prove that $\vert v_k^\eps\vert_{E^{s}_{(1)}}$ is bounded from above by
the r.h.s. of the estimate given in the statement of the lemma. The lemma
thus follows by finite induction: we just have to prove that the
desired estimate on
$v_k^\eps$ holds in $E^{s}_{(l+1)}$ ($1\leq l\leq p$) if it holds in 
$E^{s}_{(l'+1)}$, for all $l'<l$. Since moreover 
$\vert v_k^\eps\vert_{E^s_{(l+1)}}\leq \vert v_k^\eps\vert_{E^s_{(l)}}
+\vert (\eps\dt)^l\dt v_k^\eps\vert_{X^{s-(l+1)m}_T}$, we are reduced
to prove that this latter term is bounded from above by the r.h.s. of
the estimate given in the lemma.
From the equation one gets
\begin{equation}
	\label{step3}
	\vert (\eps\dt)^l\dt v_k^\eps\vert_{X^{s-(l+1)m}_T}\leq
	\vert (\eps\dt)^l (\cGeu[t,u_k^\eps] v_k^\eps)\vert_{X^{s-(l+1)m}_T}
	+\vert (\eps\dt)^l\Phi_1(u_k^\eps)\vert_{X^{s-(l+1)m}_T}. 
\end{equation}
Proceeding exactly as for the obtention of (\ref{step2}) (with $j=1$) -- 
but replacing
$\cFe$ by $\cGe$ (which is possible thanks to Lemma \ref{lmb1}), 
$\uu$ by $u_k^\eps$ and $\uv_1$ by $v_k^\eps$ -- one gets
\begin{equation}
	\vert (\eps\dt)^l (\cGeu[t,u_k^\eps] v_k^\eps)\vert_{X^{s-(l+1)m}_T}
	\leq C(s,T,\vert u_k^\eps\vert_{X^{s_0+(l+1)m}_{(l)}})
	\label{step4}
	\big(
	\vert v_k^\eps\vert_{X^{s}_{(l)}}
	+\vert u_k^\eps \vert_{X^{s}_{(l)}}
	\vert v_k^\eps\vert_{X^{s_0+(l+1)m}_{(l)}}\big).
\end{equation}
It follows therefore from (\ref{step3}) and (\ref{step4}) that
$\vert (\eps\dt)^l\dt v_k^\eps\vert_{X^{s-(l+1)m}_T}$ is bounded from 
above by
$$
	C(\eps_0,s,T,\vert u_k^\eps\vert_{E^{s_0+(l+1)m}_{(l)}})
	\big(\vert v_k^\eps\vert_{E^{s}_{(l)}}
	+\vert u_k^\eps \vert_{E^{s}_{(l)}}
	\vert v_k^\eps\vert_{E^{s_0+(l+1)m}_{(l)}}\big)
	+\vert \Phi_1(u_k^\eps)\vert_{X^{s-m}_{(l)}}.
$$
and using the induction property thus gives the result.
\end{proof}

Similarly, Lemmas \ref{prooflem2}-\ref{prooflem4} must be replaced
by the following generalizations to the case $p>0$; for Lemma \ref{prooflem2},
this is done in the following lemma.
\begin{lemmabis}\label{lmb3}
	Let $D\geq m +d_1'$ and $s\geq s_0+m$.
	If, for some $M>0$, one has 
	$\vert u_j^\eps\vert_{E^{s+D}_{(p+1)}}\leq M$ ($j=k,k+1$),
	then
	$$
	\vert \Phi^\eps(u_{k+1}^\eps)\vert_{F^{s+d_1'}_{(p)}}\leq 
	C(\eps_0,s,T,M)\big( \theta_k^{m +d_1'-D}
	+\vert v_k^\eps\vert_{E^{s+D}_{(p+1)}}\big)
	\vert v_k^\eps\vert_{E^{s+D}_{(p+1)}}.
	$$ 	
\end{lemmabis}
\begin{proof}
One must add to the proof of Lemma \ref{prooflem2} a control of 
$\vert (\eps\dt)^i \phi_1(u_{k+1}^\eps)\vert_{X_T^{s+d_1'-im}}$ 
($i\leq p$).
Owing to (\ref{proofl1}), we are reduced to control 
$\eps^i\dt^i E_1$ and 
$\eps^i\dt^i E_2$ in $X_T^{s+d_1'-im}$. From the explicit expression
of $E_1$ given in (\ref{eqE1}) and since Lemma \ref{lmb1}' allows
us to use Assumption \ref{assb2}' with $\cFe$
replaced by $\cGe$, one gets 
$$
	\vert (\eps\dt)^i E_1\vert_{X_T^{s+d_1'-im}}
	\leq 
	C(\eps_0,s,T,M) \vert v_k\vert^2_{E^{s+D}_{(i)}};
$$
similarly, one gets from (\ref{eqE2}) that
$$
	\vert (\eps\dt)^i E_2\vert_{X^{s+d_1'-im}_T}
	\leq C(\eps_0,s,T,M)\theta_k^{m +d_1'-D}\vert v_k^\eps\vert_{E^{s+D}_{(i+1)}},
$$
and the lemma follows.
\end{proof}
The generalization of Lemma \ref{prooflem3} is straightforward thanks to
Lemma \ref{existbis}':
\begin{lemmabis}\label{lmb4}
	Let $D\geq d_1$ and $s\geq s_0+(p+1)m$.
	If, for some $M>0$, one has 
	$\vert u_{k+1}^\eps\vert_{E^{s+D}_{(p+1)}}\leq M$,
	then
	$$
	\vert v_{k+1}^\eps\vert_{E^{s}_{(p+1)}}\leq
	C(\eps_0,s,T,M) \vert \Phi^\eps(u_{k+1}^\eps)\vert_{F^{s+d_1'}_{(p)}}.
	$$ 	
\end{lemmabis}
Finally, Lemma \ref{prooflem4} is generalized as follows:
\begin{lemmabis}\label{lmb5}
	Let $P\geq D\geq \delta$ and
	$s\geq s_0+(p+1)m$.
	If, for some $M>0$, one has 
	$\vert u_{k}^\eps\vert_{E^{s+D}_{(p+1)}}\leq M$ and
	$\vert (h^\eps,\uinit)\vert_{F^{s+P-m}_{(p)}}\leq M$, 
	then,  
	$$
	\vert u_{k+1}^\eps\vert_{E^{s+P}_{(p+1)}}\leq
	C(\eps_0,s,T,M)(1+\theta_k^\delta)
	(1
	+\vert u_k^\eps\vert_{E^{s+P}_{(p+1)}}).
	$$
	If moreover $\vert u_{k+1}^\eps\vert_{E^{s+D}_{(p+1)}}\leq M$, 
	then one also
	has
	$$
	\vert v_{k+1}^\eps\vert_{E^{s+P-\delta}_{(p+1)}}
	\leq C(\eps_0,s,T,M)(1
	+\vert u_{k+1}^\eps\vert_{{E^{s+P}_{(p+1)}}}).
	$$
\end{lemmabis}
\begin{proof}
Thanks to Lemma \ref{existbis}' one can generalize (\ref{proofl8})
for all $r\geq s_0+pm$ as
$$
	\vert \Phi^\eps(u)\vert_{F^r_{(p)}}\leq
	C(\vert u\vert_{E^{s_0+m }_{(p+1)}})\vert u\vert_{E^{r+m }_{(p+1)}}
	+\cst\vert (h^\eps,\uinit)\vert_{F^r_{(p)}}
$$
while (\ref{proofl9}) can be straightforwardly replaced by
$$
	\vert u_{k+1}^\eps\vert_{E^{s+P}_{(p+1)}}\leq
	\vert u_{k}^\eps\vert_{E^{s+P}_{(p+1)}}
	+\cst \theta_k^{\delta}
	\vert v_k^\eps\vert_{E^{s+P-\delta}_{(p+1)}}.
$$
Using Lemma \ref{existbis}' instead of Lemma \ref{exist}, one concludes 
as in the proof of Lemma \ref{prooflem4}.
\end{proof}
The rest of the proof of the theorem is similar to the
proof of Theorem \ref{th1}.
\end{proof}

\subsection{A few remarks}\label{sectremarks}

\subsubsection{Dependence on other parameters}\label{other}

The mappings $\cLe$ and $\cFe$ which appear in the IVP (\ref{intro1})
may also depend on other parameters than $\eps$. 
Theorems \ref{th1} (or \ref{th1b}')
still hold, with an existence time independent of all these
parameters as soon as all the constants which appear
in Assumptions \ref{assintro1}-\ref{assintro3} (or \ref{assb1}'-\ref{assb3}')
are uniform with respect
to these parameters (see Remark \ref{remstab} below for such an example).

\subsubsection{Restricting the range of the assumptions}\label{restr}

It sometimes occurs that Assumptions \ref{assb2}' (resp. \ref{assb3}' )
does not hold for
all $u\in X^{s+m}$  (resp. $u\in X^{s+d_1}_{(p+1)}$) 
but only for $u\in \Omega_0$, with $\Omega_0$ an open subset of 
$X^{s+D}$ (resp. $X^{s+D}_{(p+1)}$). If for all $\theta_0$ one can find
$u_0^\eps\in \Omega_0$ such that conditions (i)$_0$, (ii)$_0$ and (iii)$_0$ of
the induction proof of Theorems \ref{th1} and \ref{th1b}' are
satisfied, then these theorems remain true. Indeed, by choosing $\theta_0$
large enough, one can make $\vert v_k^\eps\vert_{E^{s+D}_{(p+1)}}$
($k\in\N$)
small enough to have $u_k^\eps:=u_0^\eps+\sum_{l=0}^{k-1}S_lv_l^\eps\in 
\Omega_0$.\\
In particular, \emph{the theorems still hold if the
$u_0^\eps$ provided by Lemma \ref{lmmil} belongs to $\Omega_0$.}
\begin{example}
	For the Serre and Green-Naghdi equations below, such restrictions
	on the range of validity of Assumptions \ref{assintro2} and 
	\ref{assintro3} are imposed by the ``nonzero depth condition''
	(\ref{nonzero}). The comment above shows that these
	restrictions are without consequence provided that
	(\ref{nonzero}) is initially satisfied.
\end{example}

\subsubsection{Approximate linearization}\label{sectappr}

The linear initial value problem (\ref{assIVP}) considered in 
Assumption \ref{assintro3} is the exact linearization of
(\ref{intro1}). One could replace it by an approximate linearization
in the following sense (using the same notations as in Theorem \ref{th1},
and with ${\mathcal R}[t,u]:=\dt u+\frac{1}{\eps}\cLe(t)u+\cFe[t,u]-h^\eps$):
\begin{proposition}\label{propap}
	Let ${\bf L}\in C([0,T]\times X^{s+D}; {\mathfrak L}(X^{s+m},X^s))$
	($s\geq s_0+m$) be such that 
	for all $t\in[0,T]$ and $u,v\in X^{s+D}$, one has
	$$
	\big\vert \cFeu[t,u]v-{\bf L}[t,u]v\big\vert_{X^{s+d_1'}_T}
	\leq C(s,T,\vert u\vert_{s+D})
	\vert {\mathcal R}[t,u]\vert_{X^{s+D}_T}\vert v\vert_{X^{s+D}_T}.
	$$
	Then Theorem \ref{th1} still holds if the IVP (\ref{assIVP})
	is replaced in Assumption \ref{assintro3} by
	$$
	\dsp \dt \ve +\frac{1}{\eps}\cLe(t) \ve+{\bf L}[t,\ue]\ve=f^\eps,
	\qquad
	\ve\init=g^\eps.
	$$
\end{proposition}
\begin{remark}
	In other words, the proposition states that
	one can replace the derivative of $\cFe$ in Assumption 
	\ref{assintro3}  by another linear operator, provided that
	the difference between both operators vanishes on the
	set of solutions of (\ref{intro1}). This trick can sometimes
	simplify the computations (see for instance \cite{IOP} and
	Remark \ref{remvanish} below).
\end{remark}
\begin{proof}
In the proof of Theorem \ref{th1}, we only use once the fact that
the IVP (\ref{assIVP}) is the exact linearization of (\ref{intro1}):
in the derivation of (\ref{eqE1}) (otherwise, we only
use the estimate provided by Assumption \ref{assintro1}). 
Replacing $\cFeu[t,\ue]$
by ${\bf L}[t,\ue]$ in Assumption \ref{assintro3} thus implies the
following modification of (\ref{proofl1}):
$$
	\Phi_1(u^\eps_{k+1})=E_1+E_2+E_3,
$$
where $E_1$ and $E_2$ are the same as in (\ref{proofl1}), while 
$E_3$ is given by
$$
	E_3:=\cGeu [\cdot,u_k^\eps]v_k^\eps-
	{\bf G}[\cdot,u_k^\eps]v_k^\eps,
	\quad\mbox{ with }\quad
	{\bf G}[\cdot,u_k^\eps]:=U^\eps(-t){\bf L}[\cdot,U^\eps(t) u_k^\eps].
$$
From Assumption \ref{assintro1} and the assumption made on ${\bf L}$,
one gets
$$
	\vert E_3\vert_{X_T^{s+d_1'}}\leq
	C(s,T,M)\vert \Phi_1(u_k^\eps)\vert_{X_T^{s+D}}
	\vert v_k^\eps\vert_{X_T^{s+D}},
$$
where $M$ is as in Lemma \ref{prooflem2}. Therefore, the result given in Lemma
\ref{prooflem2} remains valid provided that 
$\vert \Phi_1(u_k^\eps)\vert_{F^{s+D}}\leq C(s,T,M)\theta_k^{-q}$,
with $q=D-m-d_1'$. This point can be added without difficulty
to the induction proof of Theorem \ref{th1}:
\begin{itemize}
	\item The property is true for $k=0$ if $T$ is small
	enough (Lemma \ref{lmmil});
	\item If the property is true for some $k\in\N$, then
	we just saw that Lemma \ref{prooflem2} remains true, so 
	that $\vert \Phi_1(u_k^\eps)\vert_{X_T^{s+d_1'}}
	\leq C(s,T,M)\theta_k^{-2q}=C(s,T,m)\theta_{k+1}^{-2q/r}$;
	\item We also get that $\vert \Phi_1(u_{k+1}^\eps)\vert_{X_T^{s+P-m}}
	\leq C(s,T,M)\theta_{k+1}^{\alpha}$ (from (i)$_{k+1}$);
	\item The estimate $\vert \Phi_1(u_{k+1}^\eps)
	\vert_{X_T^{s+D}}\leq C(s,T,M)\theta_{k+1}^{-q}$ is then recovered 
	by interpolation between the two estimates above.
\end{itemize}
It follows therefore that the proof is not affected by replacing
$\cFeu$ by its approximation ${\bf L}$ in (\ref{assIVP}), which
proves the proposition.
\end{proof}

\subsection{A stability property}\label{sectstab}

We prove here a stability property for the IVP (\ref{intro1}) which is
very useful for the justification of asymptotic approximations of
the exact solution. More precisely, assume that there exists
an approximate solution $\uapp$ to (\ref{intro1}) in the sense that
\begin{equation}
	\label{intro1app}
	\left\lbrace
	\begin{array}{l}
	\dsp \dt \uapp +\frac{1}{\eps}\cLe(t) \uapp+\cFe[t,\uapp]=h^\eps
	+\iota_\eps R^\eps\\
	\uapp\,\init=\uinit+\iota_\eps r^\eps,
	\end{array}\right.
\end{equation}
with $\iota_\eps>0$ and  $(R^\eps,r^\eps)_\Rge$ bounded in some
appropriate space. Our goal here is to prove that there exists
an exact solution $\uue$ to (\ref{intro1}) and that the error made
by the approximation, namely $\uue-\uapp$, remains ``small''.  
An application of the following corollary is given in Theorem \ref{justif}
below.
\begin{corollary}\label{cor1}
	Let the assumptions of Theorem \ref{th1b}' be satisfied and
	$s\geq s_0+(p+1)m$. Assume moreover that $(\uapp)_\Rge$
	and $(R^\eps,r^\eps)_\Rge$ are  bounded in 
	$X^{s+P}_{(p)}$ and $F^{s+P}_{(p)}$ respectively.\\
	There exist $0<\underline{T}\leq T$
	and a unique family 
	$(\underline{u}^\eps)_{\Rge}$  bounded in 
	$C([0,\underline{T}];X^{s+D})$ 
	and solving the IVPs (\ref{intro1})$_\Rge$. Moreover, one has
	$$
	\vert \uue-\uapp\vert_{X^{s+D}_T}\leq \cst \iota_\eps,
	$$
	and one can take $\underline{T}=T$ if $\iota_\eps$ is small enough.
\end{corollary}
\begin{proof}
Let us seek an exact solution $\ue$ under the form 
$\ue=\uapp+\iota_\eps e^\eps$, which is equivalent to solving the IVP
$$
	\left\lbrace
	\begin{array}{l}
	\dt e^\eps+\frac{1}{\eps}\cLe(t)e^\eps
	+\ucFe[t,e^\eps]
	=-R^\eps\\
	e^\eps\init=-r^\eps,
	\end{array}\right.
$$
with
$$
	\ucFe[t,u]:=\iota_\eps^{-1}
	\big( \cFe[t,\uapp+\iota_\eps u]-\cFe[t,\uapp]\big).
$$
\begin{lemma}
	The mapping $\ucFe$ satisfies Assumption \ref{assb2}' for
	all $s$ such that the family
	$(\vert \uapp\vert_{X^{s+(i+1)m}_{(p)}})_\Rge$ 
	is  bounded.
\end{lemma}
\begin{proof}
Let us first prove that $\ucFe$ satisfies the estimates given in 
Assumption \ref{assb2}' when $j=0$.
For all $0\leq i\leq p$, one computes
that $\eps^i\underline{\mathcal F}^{\eps(i)}[t,u]$ is a sum of terms of 
the form
$$
	\iota_\eps^{-1} \big(\eps^k{\mathcal F}^{\eps(k)}_{(j')}
	[t,\uapp+\iota_\eps u]
	-\eps^k{\mathcal F}^{\eps(k)}_{(j')}[t,\uapp]\big)
	([(\eps\dt)\uapp]^{\alpha_1},\dots,
	[(\eps\dt)^i\uapp]^{\alpha_{i}}),
$$
with $k+\alpha_1+\dots+i\alpha_i=i$ and $j'=\alpha_1+\dots+\alpha_i$, 
so that we can use 
Assumption \ref{assb2}' to get
$$
	\supT \vert \eps^i\underline{\mathcal F}^{\eps(i)}[t,u]\vert_{s-im}
	\leq C(s,t,\vert \uapp\vert_{X_{(p)}^{s+(i+1)m}})\vert u\vert_{s+m},
$$
which proves the case $j=0$ since we assumed that
$(\vert \uapp\vert_{X_{(p)}^{s+(i+1)m}})_\Rge$ is  bounded.\\
Since for all $j\geq 1$ one has 
$$
	\underline{\mathcal F}^\eps_{(j)}[t,u](v_1,\dots,v_j)=
	\iota_\eps^{(j-1)}{\mathcal F}^\eps_{(j)}[t,\uapp+\iota_\eps u]
	(v_1,\dots,v_j),
$$
the case $j\geq 1$ of the Assumption  follows easily.
\end{proof}
Thanks to the lemma and to the assumptions made in the statement of the
corollary, one can use Theorem \ref{th1b}' with $\cFe$ replaced by
$\ucFe$, and the first part of the corollary is proved.\\
We now prove that it is possible to take
$T=\underline{T}$ when $\iota_\eps$ is small enough. 
Instead of taking the first iterate $u_0^\eps$ of
the converging sequence $(u_k^\eps)_k$ as given by Lemma \ref{lmmil}, we can
take $u_0^\eps=\uapp$. Instead of shrinking $T$ to prove the first step
of the induction as in the proof of Theorem \ref{th1}, one must restrict
to $\iota_\eps$ small enough. This shows that an exact 
solution to (\ref{intro1})
exists over $[0,T]$. In order to prove the error estimate, one 
proceeds as for the uniqueness part of 
Theorem \ref{th1}.
\end{proof}

\section{Application to the Green-Naghdi and Serre equations}\label{sectGN}

This section is devoted to the proof of a well-posedness and stability 
result for the Green-Naghdi and Serre equations, which are among the most
commonly used models in coastal oceanography.

\subsection{The equations}\label{secteq}

The Green-Naghdi and Serre equations describe the motion of a layer of 
incompressible and irrotational fluid under the influence of gravity and 
under some assumptions on the physical regime considered. Defining the
dimensionless parameters $\mu$ and $\eps$ as
$$
	\sqrt{\mu}:=\frac{\mbox{mean depth}}{\mbox{typical wave-length}}
	\quad\mbox{ and }\quad
	\eps:=\frac{\mbox{surface and bottom variations}}
	{\mbox{mean depth}},
$$
the Green-Naghdi and Serre regimes can be characterized as follows:
\begin{itemize}
\item Green-Naghdi regime: $\mu\ll 1$ and  $\eps\sim 1$;
\item Serre regime: $\mu\ll 1$ and $\eps\sim \sqrt{\mu}$.
\end{itemize}
A rigorous derivation of the Green-Naghdi and Serre models 
is performed in \cite{AlvarezLannes}, to which we refer for more details. 
In nondimensionalized variables, the surface is parameterized at time $t$
by $\zeta(t,X)$ ($X\in\R^2$), while the bottom is parameterized by $b(X)$.
Denoting by $V(t,X)\in \R^2$ the vertically averaged horizontal component
of the velocity field at time $t$, the equations read (with $\eps=1$ for
the Green-Naghdi equations and $\eps=\sqrt{\mu}$ for the Serre equations):
\begin{equation}\label{GN0}
	\left\lbrace
	\begin{array}{l}
	(h+\mu{\mathcal T}[h,\eps b])\dt V+h\nabla\zeta
	+h\eps(V\cdot\nabla)V\\
	\qquad \qquad+\mu\eps \Big[
	\frac{1}{3}\nabla\big(h^3{\mathcal D}_{V}\mbox{div}(V)\big)+
	\cQ[h,\eps  b](V)\Big]=0	\\
	\dt\zeta+\nabla\cdot(hV)=0,
	\end{array}\right.
\end{equation}
where $h:=1+\eps(\zeta-b)$ while the linear operators
${\mathcal T}[h,b]$ and ${\mathcal D}_{V}$ and the quadratic form 
$\cQ[h,b](\cdot)$ are 
defined as
\begin{eqnarray*}
	{\mathcal T}[h,b] V&:=&
	-\frac{1}{3}\nabla(h^3\nabla\cdot V)
	+\frac{1}{2}\big[
	\nabla(h^2\nabla b \cdot V)-h^2\nabla b \nabla\cdot V\big]
	+h\nabla b\nabla b\cdot V,\\
		{\mathcal D}_V&:=&-(V\cdot \nabla)+\mbox{div}(V)\\
	\cQ[h,b](V)&:=&	\frac{1}{2}\nabla\big(h^2(V\cdot\nabla)^2b\big)+
	h\big(\frac{h}{2}{\mathcal D}_{V}\mbox{div}(V)+(V\cdot\nabla)^2b\big)\nabla b.
\end{eqnarray*}

\subsection{Well-posedness of the Serre and Green-Naghdi equations}\label{sectWP}

Under the ``nonzero depth condition''
\begin{equation}
	\label{nonzero}
	\exists h_0>0,\qquad \inf_{\R^2}h\geq h_0,
	\qquad (h=1+\eps(\zeta-b))
\end{equation}
and after defining the spaces
\begin{equation}\label{defX}
	X^s:=\{(V,\zeta)\in H^s(\R^2)^2\times H^s(\R^2),\mbox{ such that }
	\vert \nabla\cdot V\vert_{H^s}<\infty\},
\end{equation}
endowed with the norm
\begin{equation}\label{defnorm}
	\vert (V,\zeta)\vert_{X^s}:=\Vert V\Vert_s+\vert \zeta\vert_{H^s},
	\quad\mbox{ with }\quad \Vert V\Vert_{s}:=\vert V\vert_{H^s}
	+\sqrt{\mu} \vert \nabla\cdot V\vert_{H^s},
\end{equation}
one can prove the following well-posedness result on the Serre 
($\eps=\sqrt{\mu}$) and Green-Naghdi ($\eps=1$) equations:
\begin{theorem}[Well-Posedness of the Serre and Green-Naghdi equations]
	\label{WP}
 	Let $t_0>1$, 
	$\eps=\sqrt{\mu}$ (Serre) or $\eps=1$
	(Green-Naghdi), and $s\geq t_0+2$.\\
	Let also $(V_0^\mu,\zeta_0^\mu)_{0<\mu<1}$ be
	bounded in $X^{s+38}$ and
	satisfy (\ref{nonzero}).\\
	Then there
	exists $T>0$ such that the Serre or
	Green-Naghdi  equations (\ref{GN0}) admit a unique
	family of solutions $(V^\mu,\zeta^\mu)_{0<\mu<1}$
	bounded in $C([0,\frac{T}{\eps}];X^{s+4})$ and
	with initial condition $(V_0,\zeta_0)_{0<\mu<1}$.
\end{theorem}
\begin{remark}
The spaces $X^{s+38}$ and $X^{s+4}$ correspond to the
spaces $X^{s+P}$ and $X^{s+D}$ of Theorem \ref{th1} (one can 
check that $D=4$, $P=38$ is an admissible 
choice when $m=d_1=2$ and $d_1'=0$).
\end{remark}
\begin{remark}
For
$1D$ surface waves, flat bottoms ($b=0$) and in the Green-Naghdi
scaling ($\eps=1$), Y. A. Li \cite{LiCPAM} uses precise 
estimates on the inverse
operator $(h+\mu\cT[h,\eps b])^{-1}$ to
obtain a well-posedness result in $H^{s+1}(\R)^2\times H^s(\R)$, 
with $s>3/2$ by a standard fixed point technique.
 It is not clear whether these
techniques can be adapted to the $2DH$ case: the identity
$X^s=H^{s+1}(\R^d)^d\times H^s(\R^d)$ is false when $d=2$, and
the smoothing properties of $(h+\mu\cT[h,\eps b])^{-1}$ can only
be used to control the derivatives of $V$ which are in divergence form.
\end{remark}
\begin{proof} We only prove the theorem in the Serre scaling 
($\eps=\sqrt{\mu}$), which is the most difficult one because the
existence time provided by the theorem is ``large'' (of order $O(1/\eps)$).
The modifications to prove the theorem in the Green-Naghdi scaling
are straightforward.\\
\emph{For the sake of simplicity, we write $\oper$ instead of $h+\mu\cT[h,\eps b]$ 
when no confusion is possible}.
It can be remarked that the operator $\oper$ is $L^2$ self-adjoint; since
moreover, one has
\begin{eqnarray*}
	\lefteqn{ (h+\mu\cT[h,\eps b])V,V)=(hV,V)}\\
	&+&\mu\big(h(\frac{h}{\sqrt{3}}\nabla\cdot V-\frac{\sqrt{3}}{2}\nabla b\cdot V),\frac{h}{\sqrt{3}}\nabla\cdot V-\frac{\sqrt{3}}{2}\nabla b\cdot V\big)+\frac{\mu}{4}(h\nabla b\cdot V,\nabla b\cdot V),
\end{eqnarray*}
and using the assumption that $ \inf_{\R^2}h\geq h_0$, one deduces that
\begin{equation}
	\label{GN1}
	(\oper V,V)\geq E[\eps b](V)^2,
\end{equation}
with $E[b](V)^2:=h_0\vert V\vert_{L^2}^2
+\mu h_0\big\vert \frac{h}{\sqrt{3}}\nabla\cdot V
-\frac{\sqrt{3}}{2}\nabla b\cdot V\big\vert_{L^2}^2
+\mu\frac{h_0}{4}\big\vert \nabla b\cdot V\big\vert_{L^2}^2$.

It follows that $\oper$ has a self-adjoint, 
positive, inverse bounded on $L^2(\R^2)^2$ and the equations (\ref{GN0})
can be recast under the form
\begin{equation}\label{GN3}
	\dt u+\cL u+\eps\cFe[u]=0,
\end{equation}
with $u=(V,\zeta)^T$, $\cL=\left(\begin{matrix}0 & \nabla\\ \dive & 0\end{matrix}\right)$, $\cFe[\cdot]=({\mathcal F}^\eps_1[\cdot],{\mathcal F}^\eps_2[\cdot])^T$ and
\begin{eqnarray*}
	{\mathcal F}_1^\eps[u]&=&\frac{1}{\eps}\big(\oper^{-1}h-1\big)\nabla \zeta
	+\oper^{-1}h(V\cdot\nabla)V\\
	&+&\mu\oper^{-1}
	\Big[
	\frac{1}{3}\nabla\big(h^3{\mathcal D}_{V}\mbox{div}(V)\big)+
	\cQ[h,\eps  b](V)\Big],\\
	{\mathcal F}_2^\eps[u]&=&\nabla\cdot\big((\zeta-b)V\big).
\end{eqnarray*}
The existence result stated by the theorem gives a time interval 
$[0,\frac{\uT}{\eps}]$ for (\ref{GN3}). 
Rescaling time as $t\leadsto t/\eps$, this
is equivalent to solve the following equation on the time interval
$[0,\uT]$:
\begin{equation}\label{GN4}
	\dt u+\frac{1}{\eps}\cL u+\cFe[u]=0;
\end{equation}
this latter formulation is of the form (\ref{intro1}), and the result
thus follows from Theorem \ref{th1}, provided that 
Assumptions \ref{assintro1}-\ref{assintro3} are satisfied. The rest
of the proof is devoted to check that these assumptions are satisfied
in the Banach scale $X^s$ defined in (\ref{defX}), with $s_0=t_0$,
$m=d_1=2$, and $d_1'=0$.
\begin{remark}\label{remstab}
	In the Green-Naghdi scaling ($\eps=1$), the parameter 	
	$\mu$ cannot be expressed in terms of $\eps$. As seen
	in Section \ref{other}, in order for the theorem to be valid
	uniformly with respect to $\mu\in(0,1)$, we must check that
	Assumptions \ref{assintro1}-\ref{assintro3} are satisfied
	uniformly in $\mu\in(0,1)$.
\end{remark}
\begin{remark}
	We always assume implicitly that the nonzero depth condition
	(\ref{nonzero}) is satisfied. As explained in Section \ref{restr},
	this is implied by the assumption that (\ref{nonzero})
	is satisfied at $t=0$.
\end{remark}

It follows immediately from the definition of 
$\cL$ (independent of time here) that Assumption \ref{assintro1} is satisfied. In order to check the other assumptions, we need some preliminary results.\\ 
The next lemma gathers some general estimates; the first one is a classical
Moser tame product estimate, the second one is a generalized Kato-Ponce
commutator estimate (note that the estimate depends on $f$ only
through its gradient, see Ths. 3 and 6 of \cite{LannesJFA}), and the last
one is a classical ``quasilinear type'' estimate 
(e.g. Chapter II.C of \cite{AlinhacGerard}).
\begin{lemma}
	\label{propprel1}
	Let $t_0>d/2$ and $s\geq 0$.\\
	{\bf i.} For all
	$f,g\in H^{s}\cap H^{t_0}(\R^d)$, one has, 
	using notation (\ref{nota3}),
	$$
	\vert fg\vert_{H^s}\lesssim
	\vert f\vert_{H^{t_0}}\vert g\vert_{H^s}
	+ \langle \vert f\vert_{H^{s}}
	\vert g\vert_{H^{t_0}}\rangle_{s> t_0};
	$$
	{\bf ii. }Let $r\in\R$ be such that $-t_0<r\leq t_0+1$.
	Then, for all $f\in H^{t_0+1}\cap H^{s+r}(\R^d)$
	and $u\in H^{t_0}\cap H^{s+r-1}(\R^d)$,
	$$
	\big\vert
	[\Lambda^s,f]u\big\vert_{H^r}
	\lesssim \vert \nabla f\vert_{H^{t_0}}\vert u\vert_{H^{s+r-1}}
	+\left\langle  \vert \nabla f\vert_{H^{s+r-1}}\vert u\vert_{H^{t_0}}
	\right\rangle_{s> t_0+1-r}.
	$$
	{\bf iii.} Let $N\in\N$, and ${\bf P}$ be a first order
	differential operator on $L^2(\R^d)^N$ with anti-adjoint
	principal part: ${\bf P}:=\sum_{j=1}^d P_j(x)\partial_j+P_0(x)$,
	with $P_j(x)$ symmetric for $j=1,\dots, d$ and $x\in\R^d$. Then, 
	for all $U\in H^s\cap H^{t_0+1}(\R^d)^N$,
	$$
	\vert \big(\Lambda^s{\bf P} U,\Lambda^s U\big)\vert\lesssim
	\Big(\Vert {\bf P}\Vert_{H^{t_0+1}}\vert U\vert_{H^s}
	+\big\langle \Vert {\bf P}\Vert_{H^s}\vert U\vert_{H^{t_0+1}}
	\big\rangle_{s>t_0+1}\Big)\vert U\vert_{H^s},
	$$
	with $\Vert {\bf P}\Vert_{H^s}:=\sum_{j=0}^d\vert P_j\vert_{H^s}$.
\end{lemma}
The following lemma gives some properties on $\imfT$ which are necessary
to check Assumption \ref{assintro2} (recall that the norm
$\Vert\cdot\Vert_s$ has been defined in (\ref{defnorm})).
\begin{lemma}\label{lmcheck}
	The following estimates hold for all $s\geq 0$ and uniformly with 
	respect to $\mu\in (0,1)$ (and $\eps=1$ or $\eps=\sqrt{\mu}$):
	\begin{enumerate}
	\item[i. ] $\displaystyle \Vert \imfT V\Vert_s
	\leq
	c_0
	\big(\vert V\vert_{H^s}
	+\big\langle (\vert h\vert_{H^s}
	+\vert \nabla b\vert_{H^s})\vert V\vert_{H^{t_0}}
	\big\rangle_{s>t_0+1}\big)
	$;
	\item[ii. ]
	$\displaystyle
	\sqrt{\mu}\big\vert \oper^{-1}\nabla \zeta \vert_{H^s}
	\leq
	c_0
	\big(\vert \zeta \vert_{H^s}
	+\big\langle (\vert h\vert_{H^s}
	+\vert \nabla b\vert_{H^s})\vert \zeta \vert_{H^{t_0+1}}
	\big\rangle_{s>t_0+1}\big)
	$;
	\item[iii. ]
	$
	\displaystyle 
	\frac{1}{\eps}\big\vert \big(\oper^{-1}h-1)V\big\vert_{H^s}
	\leq
	c_0
	\big(\vert V\vert_{H^{s+1}}
	+\big\langle (\vert h\vert_{H^{s+1}}
	+\vert \nabla b\vert_{H^{s+1}})\vert V\vert_{H^{t_0+1}}
	\big\rangle_{s>t_0}\big)
	$,
	\end{enumerate}
	where $c_0$ is a constant depending only on 
	$\frac{1}{h_0}$, $\vert h\vert_{H^{t_0+1}}$ and 
	$\vert\nabla b\vert_{H^{t_0+1}}$.
\end{lemma}
\begin{proof}
{\bf i.} Remark first that 
\begin{equation}\label{rel}
	\vert V\vert_{L^2}^2\leq \frac{1}{{h_0}}E[\eps b](V)^2
	\quad \mbox{ and }\quad
	{\mu}\vert\nabla\cdot V\vert_{L^2}^2\leq 
	C(\frac{1}{h_0},\vert \nabla b\vert_{L^\infty})E[\eps b](V)^2;
\end{equation}
replacing $V$ by $\imfT V$ in the above expressions and using
(\ref{GN1}) shows that
$$
	\vert \imfT V\vert_{L^2}^2\leq \frac{1}{h_0}
	(V,\imfT V)
	\quad\mbox{ and }\quad
	{\mu}\vert \nabla\cdot \imfT V\vert_{L^2}^2
	\leq C(\frac{1}{h_0},\vert \nabla b\vert_{L^\infty})
		(V,\imfT V).
$$
A simple Cauchy-Schwartz inequality thus yields
\begin{equation}\label{estnorm}
	\left\lbrace
	\begin{array}{l}
	\Vert \imfT \Vert_{L^2(\R^2)^2\to L^2(\R^2)^2}\leq 
	\frac{1}{h_0}\\
	\sqrt{\mu}\Vert \nabla\cdot \imfT \Vert_{L^2(\R^2)^2\to L^2(\R^2)}
	\leq 
	C(\frac{1}{h_0},\vert \nabla b\vert_{L^\infty}).
	\end{array}\right.
\end{equation}
Using the fact that $\mfT^{-1}$
is self-adjoint, one has
$\big(\imfT\nabla\big)^*=-\nabla\cdot\imfT$, and thus
\begin{equation}\label{estnorm1}
	\sqrt{\mu}\Vert \imfT\nabla\Vert_{L^2(\R^2)\to L^2(\R^2)^2}
	\leq 
	C(\frac{1}{h_0},\vert \nabla b\vert_{L^\infty}).
\end{equation}
We can now prove the following inequality
\begin{equation}
	E[\eps b](\Lambda^s \imfT V)\leq
	c_0
	\label{estnorm2}
	\big(\vert V\vert_{H^s}+\big\langle (\vert h\vert_{H^s}+\vert \nabla b\vert_{H^s})\vert V\vert_{H^{t_0}}\big\rangle_{s>t_0+1}\big),
\end{equation}
which, together with (\ref{rel}), obviously implies the first point 
of the lemma. 
From (\ref{GN1}), one gets the relation
\begin{eqnarray}
	\nonumber
	\lefteqn{E[\eps b](\Lambda^s \imfT V)^2\leq
	\big(\Lambda^s V,\Lambda^s \imfT V\big)
	+\big(\mfT[\Lambda^s,\imfT]V,\Lambda^s\imfT V\big)}\\
	\label{refcom}
	&=&\big(\Lambda^s V,\Lambda^s \imfT V\big)
	-
	\big([\Lambda^s,h]\imfT V,\Lambda^s \imfT V\big)
	-\mu\big([\Lambda^s,\cT]\imfT V,\Lambda^s \imfT V\big).
\end{eqnarray}
Replacing $\cT$ by its expression and integrating by parts, one gets
therefore
\begin{eqnarray*}
	\lefteqn{E[\eps b](\Lambda^s \imfT V)^2\leq
	\big(\Lambda^s V,\Lambda^s \imfT V\big)
	-\big([\Lambda^s,h]\imfT V,\Lambda^s \imfT V\big)}\\
	& &-\frac{1}{3}\big([\Lambda^s,h^3](\sqrt{\mu}\nabla\cdot \imfT V),
	\Lambda^s(\sqrt{\mu}\nabla\cdot \imfT V)\big)\\
	& &-\frac{\sqrt{\mu}}{2}\big([\Lambda^s,h^2\nabla b^T] \imfT V,
	\Lambda^s(\sqrt{\mu}\nabla\cdot \imfT V)\big)\\
	& &+\frac{\sqrt{\mu}}{2}\big([\Lambda^s,h^2\nabla b] 
	(\sqrt{\mu}\nabla\cdot \imfT V),
	\Lambda^s \imfT V\big)
	+\mu\big([\Lambda^s,h\nabla b\nabla b^T]\imfT V,
	\Lambda^s\imfT V\big).
\end{eqnarray*}
Applying Cauchy-Schwartz's inequality to every component of the r.h.s. of 
the above expression, and using Lemma \ref{propprel1}.ii and (\ref{GN1}),
one gets directly
\begin{eqnarray*}
	\lefteqn{E[\eps b](\Lambda^s \imfT V)
	\leq C(\frac{1}{h_0},\vert h\vert_{H^{t_0+1}},
	\vert \nabla b\vert_{H^{t_0+1}})}\\
	&\times& \Big(\vert V\vert_{H^s}+E[\eps b](\Lambda^{s-1} \imfT V)
	+\big\langle (\vert h\vert_{H^s}+\vert\nabla b\vert_{H^s})E[\eps b](\Lambda^{t_0} V)\big\rangle_{s>t_0+1}\Big),
\end{eqnarray*}
from which one deduces (\ref{estnorm2}). \\
{\bf ii.} Remark that
\begin{eqnarray}
	\nonumber
	\sqrt{\mu}\big\vert \Lambda^s\imfT\nabla\zeta\big\vert_2
	&\leq&
	\sqrt{\mu}\big\vert \imfT\nabla\Lambda^s\zeta\big\vert_2
	+\sqrt{\mu}\big\vert [\Lambda^s, \imfT]\nabla\zeta\big\vert_2\\
	\label{estnorm3}
	&\leq &C(\frac{1}{h_0},\vert \nabla b\vert_{L^\infty})
	\vert \zeta\vert_{H^s}
	+\sqrt{\mu}\big\vert [\Lambda^s, \imfT]\nabla\zeta\big\vert_2,
\end{eqnarray}
where we used (\ref{estnorm}) to control the first term of the r.h.s. For the
second term, remark that $[\Lambda^s,\imfT]=-\imfT[\Lambda^s,\mfT]\imfT$,
so that
\begin{eqnarray*}
	[\Lambda^s,\imfT]&=&-\imfT[\Lambda^s, h]\imfT
	+\frac{1}{3}(\sqrt{\mu}\imfT\nabla)[\Lambda^s,h^3]
	(\sqrt{\mu}\nabla\cdot\imfT)\\
	& &-\frac{\sqrt{\mu}}{2}(\sqrt{\mu}\imfT\nabla)[\Lambda^s,h^2\nabla b]
	\cdot\imfT+\frac{\sqrt{\mu}}{2}\imfT[\Lambda^s,\nabla b](\sqrt{\mu}\nabla\cdot\imfT)\\
	& &-\mu\imfT[\Lambda^s,\nabla b\nabla b^T]\imfT.
\end{eqnarray*}
It follows that $\vert [\Lambda^s,\imfT]\nabla\zeta \vert_{L^2}$ can  
be controlled using Lemma \ref{propprel1}.ii, the
estimates on $\Vert \imfT\Vert_{L^2\to L^2}$ and $\sqrt{\mu}\Vert \imfT \nabla\Vert_{L^2\to L^2}$ given by
(\ref{estnorm}) and (\ref{estnorm1}), and the estimate on $\Vert \imfT\Vert_{H^s\to H^s}$ and $\sqrt{\mu}\Vert \nabla\cdot\imfT\Vert_{H^s\to H^s}$ 
given by i. Together with 
(\ref{estnorm3}), this proves ii.\\
{\bf iii.} Just remark that $1-\imfT h=\imfT(\mfT-h)=\mu\imfT\cT$,
so that 
$$
	\frac{1}{\eps}\vert (1-\imfT h)V\vert_{H^s}
	\leq \frac{\sqrt{\mu}}{\eps}\vert \sqrt{\mu}\imfT \cT V\vert_{H^s}.
$$
Since $\sqrt{\mu}/\eps\leq 1$, the result follows from the first two points
and Lemma \ref{propprel1}.i.
\end{proof}

A direct consequence of Lemmas \ref{propprel1}.i and \ref{lmcheck}.i and .iii, and 
of the definition of $\cFe$, 
is that the first part of Assumption \ref{assintro2} holds with 
$s=t_0$ and $m=2$.\\
Remarking that 
$$
	d_{\zeta}(\zeta\mapsto \oper^{-1})\widetilde{\zeta}
	=-\oper^{-1} \big(\eps\widetilde{\zeta}+\mu d_{\zeta}(\zeta\mapsto \cT[h,\eps b])\widetilde{\zeta}\big)
	\oper^{-1},
$$
and using Lemmas \ref{propprel1}.i and \ref{lmcheck}.i and .ii, 
one can check that the conditions 
on the first and second derivative of $\cFe$ made in 
Assumption \ref{assintro2} are also satisfied.\\
In order to check Assumption \ref{assintro3}, we must study the
Cauchy problem associated to the linearization of (\ref{GN4}) around
some reference state $\uu:=(\uV,\uz)$:
\begin{equation}\label{GNCauchy1}
	\dt u+\frac{1}{\eps}\cL u+\cFeu[\uu]u=f,\qquad u\init=g,
\end{equation}
with $f=(F_1,f_2)^T$ and $g:=(G_1,g_2)^T$.\\
From the explicit expression
of $\cL$ and $\cFe$ given in (\ref{GN3}), and writing
$\uh:=1+\eps(\uz-b)$, $\underline{\mathcal T}:={\mathcal T}[\uh,\eps b]$,
and $\umfT:=\uh+\mu\underline{\cT}$,
one gets that (\ref{GNCauchy1})
is equivalent to
\begin{equation}\label{GNCauchy2}
	\left(\begin{array}{c} \uoper\dt V\\
	\dt \zeta \end{array}\right)
	+\left(\begin{array}{cc}
	\cN_1 & \cN_2\\
	\cN_3 &\cN_4
	       \end{array}\right) \left(\begin{array}{c} V\\ \zeta\end{array}\right)
	=
	\left(\begin{array}{c}
	\uoper F_1\\
	f_2   \end{array}\right),
\end{equation}
with initial condition $(V,\zeta)\init= (G_1,g_2)$, and
where the linear operators $\cN_j$ ($j=1,\dots,4$) are given by
\begin{eqnarray*}
	\cN_1 V&:=&\uh (\uV\cdot\nabla)V+\uh(V\cdot\nabla)\uV
	+\frac{\mu}{3}\nabla\big[\uh^3\big(\cD_{\uV}(\nabla\cdot V)
	+\cD_V(\nabla\cdot\uV)\big)\big],\\
	& &+\mu\cQ[\uh,\eps b](V,\uV),\\
	\cN_2 \zeta &:=&\frac{1}{\eps}\uh \nabla \zeta
	+\umfb\zeta
	+\mu\nabla\big(\uh\umfa \zeta  \big),\\
	\cN_3 V&:=&\frac{1}{\eps}\nabla\cdot (\uh V),\\
	\cN_4 \zeta&:=&\nabla\cdot(\zeta\uV) ,
\end{eqnarray*}
with 
\begin{eqnarray*}
	{\umfa}&:=&\eps\uh{\mathcal D}_{\uV}(\nabla\cdot\uV)+(\uV\cdot\nabla)^2b-(\nabla b-\uh\nabla)\cdot(\nabla\uz+\eps\cF^\epsilon_1[\uu]),\\ 
{\umfb}&:=& \eps(\uV\cdot\nabla)\uV-\eps\cF^\epsilon_1[\uu]
	+\mu\umfa\nabla b,
\end{eqnarray*}
and with $\cQ[\uh,\eps b](\cdot,\cdot)$ standing for the bilinear symmetric form canonically
associated to $\cQ[\uh,\eps b]$.
\begin{remark}\label{remvanish}
	Thanks to Proposition \ref{propap}, one can replace
	$\nabla\uz+\eps\cF^\eps_1[\uu]$ by $-\eps\dt\uV$ in the
	definition of $\umfa$ and $-\eps\cF^\eps_1[\uu]$ by
	$\eps\dt\uV+\nabla\uz$ in the definition of
	$\umfb$. It follows that
	$\sqrt{\mu}\vert \umfa\vert_{H^s}$ 
	and $\vert \umfb\vert_{H^s}$ are  bounded from above by
	$$
	C(\vert \uz\vert_{H^{t_0}},\Vert \uV\Vert_{t_0+1},
	\vert b\vert_{H^{t_0+1}},\Vert \eps\dt\uV\Vert_{t_0})
	\big(\vert \uz\vert_{H^{s}}+
	\Vert \uV\Vert_{s+1}+\Vert \eps\dt\uV\Vert_{s}	
	+\vert b\vert_{H^{s+1}}\big).
	$$
\end{remark}

The following lemma end the proof of the theorem since it implies that
Assumption \ref{assintro3} is also satisfied with 
	$s_0=t_0$, $m=2$, $d_1=2$ and $d_1'=0$ (recall
that the notation ${\mathcal I}^s$ is defined in (\ref{defI})):
\begin{lemma}
	 Let $T>0$, $f=(F_1,f_2)\in C([0,T];X^s)$, $g\in X^s$
	and $\uu=(\uV,\uz)\in X^{s+1}_T$ be such that 
	$\dt \uu\in X^s_T$.
	 Then for all $s\geq 0$, there exists a 
	unique solution $u=(V,\zeta)\in X^s_T$ to (\ref{GNCauchy1})
	and for all $0\leq t\leq T$,
	$$
	\vert u(t)\vert_{X^s}\leq \uc_1\big({\mathcal I}^s(t,f,g)+
	\big\langle(\vert \uu\vert_{X^{s+1}_T}+\vert \eps\dt \uu\vert_{X^s_T})
	{\mathcal I}^{t_0+1}(t,f,g)\big\rangle_{s>t_0+1}\big),
	$$
	with $\uc_1=C(T,\frac{1}{h_0},\vert \uu\vert_{X^{t_0+2}_T},
	\vert\eps \dt\uu \vert_{X^{t_0+1}_T},\vert b\vert_{H^{t_0+2}},
	\vert b\vert_{H^{s+2}})$. 
\end{lemma}
\begin{proof}
We only prove the energy estimate, since the existence/uniqueness
of a solution to the linear Cauchy problem (\ref{GNCauchy2})
can be classically
obtained by regularization techniques. Multiplying the
first equation of (\ref{GNCauchy2}) by $\umfT\Lambda^s{\uimfT}$ and the second
by $\Lambda^s$, and taking the scalar product with
$\Lambda^s V$ and $\Lambda^s \zeta$ respectively, one gets
\begin{eqnarray}
	\nonumber
	\lefteqn{\frac{1}{2}\dt (\umfT \Lambda^s V,\Lambda^s V)
	+\frac{1}{2}\dt (\Lambda^s \zeta,\Lambda^s\zeta)
	=-(\Lambda^s \cN_1 V,\Lambda^s V)-(\Lambda^s \cN_4 \zeta,\Lambda^s \zeta)}\\
	\nonumber
	& &-\big[
	(\Lambda^s \cN_2 \zeta,\Lambda^s V)
	+(\Lambda^s \cN_3 V,\Lambda^s \zeta)\big]\\
	\nonumber
	& &+([\Lambda^s,\umfT]\uimfT \cN_1 V,\Lambda^s V)
	+([\Lambda^s,\umfT]\uimfT \cN_2\zeta,\Lambda^s V)\\
	\label{devel}
	& &+(\umfT \Lambda^s F_1,\Lambda^s V)+(\Lambda^s f_2,\Lambda^s \zeta)
	+\frac{1}{2}(\dt\umfT \Lambda^s V,\Lambda^s V);
\end{eqnarray}
we now prove that (with $u=(V,\zeta)$, $\uu=(\uV,\uz)$, $f=(F_1,f_2)$,
 and $\vert\cdot\vert_{X^s}$ as defined in (\ref{defnorm}))
\begin{eqnarray}\label{estgenerale}
	\nonumber
	\lefteqn{\dt (\umfT \Lambda^s V,\Lambda^s V)
	+\dt (\Lambda^s \zeta,\Lambda^s\zeta)}\\
	&\leq& \uc_1
	\big(\vert u\vert_{X^s}+\vert f\vert_{X^s}
	+\big\langle(\vert \uu \vert_{X^{s+1}}+\vert \eps\dt\uu\vert_{X^s})
	\vert u\vert_{X^{t_0+1}}\big\rangle_{s>t_0+1}\big)
	\vert u\vert_{X^s}.
\end{eqnarray}
We thus check that all the components of the r.h.s. of (\ref{devel}) are 
bounded from above by the r.h.s. of (\ref{estgenerale}).\\
$\bullet$ Control of $(\Lambda^s\cN_1 V,\Lambda^s V)$. 
Let us first rewrite
\begin{eqnarray}\nonumber
	(\Lambda^s\cN_1 V,\Lambda^s V)&=&-\frac{1}{2}\mu I
	+(\Lambda^s {\bf A} V,\Lambda^s V)
	+\mu (\Lambda^s {\bf B} (\nabla\cdot V),\Lambda^s (\nabla\cdot V))\\
	\label{rewrite}
	& &+\mu(\Lambda^s {\bf C}V,\Lambda^s (\nabla\cdot V))
	+\mu(\Lambda^s {\bf D}(\nabla \cdot V),\Lambda^s V)
\end{eqnarray}
with 
\begin{eqnarray*}
	{\bf A}V&:=&
	\uh (\uV\cdot \nabla)\big(V+\mu(V\cdot\nabla b)\nabla b\big)
	+\uh(V\cdot\nabla)\big(\uV+\mu (\uV\cdot \nabla b)\nabla b\big)\\
	& &-\frac{\mu}{2}\uh^2(V\cdot\nabla)(\nabla\cdot\uV)\nabla b,\\
	{\bf B}&:=&
	-\frac{1}{3}\uh^3\cD_{\uV}-\frac{1}{3}\uh^3(\nabla\cdot\uV),\\
	{\bf C}V&:=& \frac{1}{3}\uh^3 (V\cdot\nabla)(\nabla\cdot\uV)
	-\frac{1}{2}\uh^2(V\cdot\nabla)(\uV\cdot \nabla b),\\
	& &+\big[\uh^2(\uV\cdot\nabla)+(\uh^2(\uV\cdot\nabla))^*\big]\nabla b\cdot V\\
	{\bf D}&:=&
	\uh^2\nabla b (\nabla\cdot \uV),
\end{eqnarray*}
and $I:=\big([\Lambda^s,\uh^2\nabla b (\uV\cdot \nabla)]\nabla\cdot V,\Lambda^s V\big)+\big([\Lambda^s,\uh^2(\uV\cdot\nabla)(\nabla b^T\cdot)]V,\Lambda^s\nabla\cdot V\big)$. One can check that $\mu I$ is bounded from above by the
r.h.s. of (\ref{estgenerale}) by applying Cauchy-Schwartz's inequality to
its two components, and then using (\ref{rel}) and Lemma \ref{propprel1}.ii.\\
Remarking also that ${\bf A}$ and ${\bf B}$ are first order differential
operators with anti-adjoint principal part, one can use Lemma \ref{propprel1}.iii and (\ref{rel}) to check that the second and third component of 
(\ref{rewrite}) are
bounded from above by the r.h.s. of (\ref{estgenerale}). Finally, we can prove 
that the same control holds on the last two components of (\ref{rewrite})
by using Cauchy-Schwartz's inequality and Lemma \ref{propprel1}.i (remark
that ${\bf C}$ and ${\bf D}$ are simple matrix and scalar valued functions).

\noindent
$\bullet$ Control of $(\Lambda^s \cN_4 \zeta,\Lambda^s \zeta)$. From the
explicit expression of $\cN_4$ and Lemma \ref{propprel1}.iii, one
obtains directly that this term is controlled by the r.h.s.
of (\ref{estgenerale}).

\noindent
$\bullet$ Control of $(\Lambda^s \cN_2 \zeta,\Lambda^s V)
+(\Lambda^s \cN_3 V,\Lambda^s \zeta)$. Integrating by parts, one gets
immediately that
\begin{eqnarray*}
	(\Lambda^s \cN_2 \zeta,\Lambda^s V)
+(\Lambda^s \cN_3 V,\Lambda^s \zeta)&=&	
	\big(\Lambda^s(\umfb\zeta),\Lambda^s V\big)
	-\big(\Lambda^s(\uh(\sqrt{\mu}\umfa)\zeta),
	\sqrt{\mu}\Lambda^s\nabla\cdot V\big)\\
	& &+\frac{1}{\eps}\big([\Lambda^s,\uh]\nabla\zeta,\Lambda^s V\big)
	+\frac{1}{\eps}\big(\nabla[\Lambda^s,\uh]V,\Lambda^s \zeta\big).
\end{eqnarray*}
The first two components can be controlled by a Cauchy-Schwartz inequality and
Lemma \ref{propprel1}.i, and the last two by Cauchy-Schwartz and 
Lemma \ref{propprel1}.ii (note that the commutator estimate provided by 
this lemma depends only on $\uh$ through $\nabla \uh=\eps\nabla (\uz-b)$
and provides therefore the $\eps$ necessary to compensate the singular
term $1/\eps$). Together with the estimates of Remark \ref{remvanish},
this shows that this term also is controlled by the r.h.s. of 
(\ref{estgenerale}).

\noindent
$\bullet$ Control of $([\Lambda^s,\umfT]\uimfT\cN_1 V,\Lambda^s V\big)$. 
Let us remark that $\cN_1 V$ can be written as
\begin{equation}
	\label{rewriteN1}
	\cN_1V=\umfT H+{\bf P}_1 V
	+{\bf P}_2(\sqrt{\mu}\nabla \cdot V),
\end{equation}
with $H:=\big[(\nabla\cdot V)\uV
-\frac{3}{2}(V\cdot\nabla b)\uV\big]$,  and where
${\bf P}_1$ and ${\bf P}_2$ are both first order differential
operators, and whose coefficients are polynomial expressions of 
$\sqrt{\mu}$, $h$, $\nabla h$, $\sqrt{\mu}\nabla\cdot\uV$,
$\nabla(\sqrt{\mu}\nabla\cdot\uV)$ and of the vectors 
$\uV$ and $\nabla b$ and their first derivatives (the exact
expression of ${\bf P}_1$ and ${\bf P}_2$ is of no importance). 
One has therefore
\begin{eqnarray*}
	[\Lambda^s,\umfT]\uimfT \cN_1 V&=&[\Lambda^s,\umfT]H
	+[\Lambda^s,\umfT]\uimfT ({\bf P}_1 V
	+{\bf P}_2(\sqrt{\mu}\nabla \cdot V))\\
	&=& [\Lambda^s,\uh]H
	+[\Lambda^s,\uh]\uimfT ({\bf P}_1 V
	+{\bf P}_2(\sqrt{\mu}\nabla \cdot V))\\	
	& &+\mu[\Lambda^s,\cT]H
	+\mu[\Lambda^s,\cT]\uimfT ({\bf P}_1 V
	+{\bf P}_2(\sqrt{\mu}\nabla \cdot V))\\
	&:=&I_1+I_2+I_3+I_4.
\end{eqnarray*}
One deduces directly from Lemmas \ref{propprel1}.ii  and \ref{lmcheck}.i
and the definition 
of $H$, ${\bf P}_1$ and ${\bf P}_2$ that
$$
	\vert I_1\vert_{L^2}+\vert I_2\vert_{L^2}\leq 
	\uc_1
	\big(\Vert V\Vert_{s}
	+\big\langle (\Vert\uV\Vert_{H^{s}}+\vert \uh\vert_{H^s}
	+\vert\nabla b\vert_{H^s})
	\Vert V\Vert_{{t_0+1}}\big\rangle_{s>t_0+1}\big),
$$
and a simple Cauchy-Schwartz inequality shows that scalar products 
$(I_1,\Lambda^s V)$ and $(I_2,\Lambda^s V)$ are controlled by
the r.h.s. of (\ref{estgenerale}). Proceeding  exactly as for the 
control of the third term of (\ref{refcom}), one can check that
the same control holds for  $(I_3,\Lambda^s V)$ and $(I_4,\Lambda^s V)$. 

\noindent
$\bullet$ Control of $([\Lambda^s,\umfT]\uimfT\cN_2 \zeta,\Lambda^s V\big)$.
Using the explicit expression of $\umfT$ and Lemma \ref{lmcheck}, and
proceeding as for the control of the third term of (\ref{refcom}), one
can bound this term from above by the r.h.s. of (\ref{estgenerale}).

\noindent
$\bullet$ Control of the last three terms of (\ref{devel}). Controlling
these terms by the r.h.s. of (\ref{estgenerale}) follows directly
from a Cauchy-Schwartz inequality (and an integration by parts for the
$O(\mu)$ component of $\umfT$ and $\dt\umfT$).

We can now conclude the proof of the lemma. Recalling that
from (\ref{GN1}), (\ref{rel}) and
the definition (\ref{defnorm}), one has
\begin{eqnarray*}
	(\umfT \Lambda^s V,\Lambda^s V)+(\Lambda^s\zeta,\Lambda^s\zeta)&\leq& 
	C(\vert \uh\vert_{L^\infty},\vert\nabla b\vert_{L^\infty})
	\vert u\vert_{X^s}^2\\
	\vert u\vert_{X^s}^2&\leq& 
	C(\frac{1}{h_0},\vert\nabla b\vert_{L^\infty})
	(\umfT \Lambda^s V,\Lambda^s V)+(\Lambda^s\zeta,\Lambda^s\zeta),
\end{eqnarray*}
one can integrate (\ref{estgenerale}) with respect to time to obtain,
for all $t\in [0,T]$,
$$
	\vert u(t)\vert_{X^s}
	\leq \uc_1 \big[\vert u(0)\vert_{X^s}+
	\int_0^t\big(\vert u\vert_{X^s}+\vert f\vert_{X^s}
	+\big\langle(\vert \uu \vert_{X^{s+1}}+\vert \eps\dt\uu\vert_{X^s})
	\vert u\vert_{X^{t_0+1}}\big\rangle_{s>t_0+1}\big)\big].
$$
Using this identity with $s=t_0+1$ shows that 
\begin{equation}\label{last}
	\vert u(t)\vert_{X^{t_0+1}}
	\leq \uc_1 \big[\vert u(0)\vert_{X^{t_0+1}}+
	\int_0^t\big(\vert u\vert_{X^{t_0+1}}+\vert f\vert_{X^{t_0+1}}\big)
	\big],
\end{equation}
and Gronwall's lemma thus yields, for all $t\in [0,T]$,
$$
	\vert u(t)\vert_{X^{t_0+1}}
	\leq \uc_1\big(\vert u(0)\vert_{X^{t_0+1}}+\int_0^t
	\sup_{0\leq t''\leq t'}\vert f(t'')\vert_{X^{t_0+1}}dt')
$$
(recall that $\uc_1$ is a generic notation whose value 
can change from one line to another); plugging this expression
into (\ref{last}) thus ends the proof of the lemma.
\end{proof}
\end{proof}

\subsection{Justification of the Serre and Green-Naghdi models}\label{sectjustif}

As said above, the Serre and Green-Naghdi are both asymptotic models
which describe the dynamics of the water-waves equations. It is not
known however whether these asymptotics are correct, in the sense that
the exact solutions to the asymptotic models provide a correct approximation
to the exact solutions of the water-waves equations. This is what we show
below: if solutions $(V^\mu_{app},\zeta^\mu_{app})$ 
to the water-waves equations exist and approximately
solve (\ref{GN0}), then their
asymptotic behavior (as $\mu\to 0$) is correctly described by
the Serre ($\eps=\sqrt{\mu}$) or Green-Naghdi ($\eps=1$) models.
\begin{theorem}[Justification of the Serre and Green-Naghdi models]
	\label{justif}
	Let $t_0>1$, 
	$\eps=\sqrt{\mu}$ (Serre) or $\eps=1$
	(Green-Naghdi), and $s\geq t_0+2$.\\
	Let also $\uT>0$ and $(V_{app}^\mu,\zeta_{app}^\mu)_{0<\mu<1}$ 
	be bounded in
	$C([0,\frac{\uT}{\eps}];X^{s+40})$ and assume that
	$h^\mu_{app}:=1+\eps(\zeta_{app}^\mu-b)$ satisfies
	(\ref{nonzero}) at $t=0$.
	If moreover
	$$
	\left\lbrace
	\begin{array}{l}
	(h^\mu_{app}+\mu{\mathcal T}[h^\mu_{app},\eps b])\dt V_{app}^\mu+h_{app}^\mu\nabla\zeta_{app}^\mu
	+\eps h_{app}^\mu(V_{app}^\mu\cdot\nabla)V_{app}^\mu\\
	\qquad \qquad+\mu\eps \Big[
	\frac{1}{3}\nabla\big({h_{app}^\mu}^3
	{\mathcal D}_{V_{app}^\mu}\mbox{\textnormal{ div}}(V_{app}^\mu)\big)+
	\cQ[h_{app}^\mu,\eps  b](V_{app}^\mu)\Big]=\mu^2 R_1^\mu	\\
	\dt\zeta_{app}^\mu+\nabla\cdot(h_{app}^\mu V_{app}^\mu)=\mu^2r_2^\mu,
	\end{array}\right.
	$$
	with  $(R_1^\mu,r_2^\mu)_{0<\mu<1}$ 
	bounded
	in $C([0,\frac{\uT}{\eps}];X^{s+38})$,
	then there exists $0<T\leq \uT$ and a unique solution 
	$(V^\mu,\zeta^\mu)_{0<\mu<1}\in C([0,T/\eps]:X^{s+4})$ 
	to (\ref{GN0}) with initial
	conditions $(V_{app}^\mu\,\init,\zeta_{app}^\mu\,\init)$.
	Moreover, one has
	$$
	\sup_{0\leq t\leq T}\big\vert (V^\mu,\zeta^\mu)	
	-(V_{app}^\mu,\zeta_{app}^\mu)
	\big\vert\lesssim \mu^2/\eps,
	$$
	uniformly with respect to $0<\mu<1$. Restricting to small enough
	values of $\mu$, one can moreover take $T=\uT$.
\end{theorem}
\begin{remark}
	The existence of a family $(V_{app}^\mu,\zeta_{app}^\mu)_{0<\mu<1}$
	having the properties assumed in the theorem is established
	in \cite{AlvarezLannes}.
\end{remark}
\begin{proof}
Since Assumptions \ref{assintro1}-\ref{assintro3} have been checked in
the proof of Theorem \ref{WP}, the result is a direct consequence
of Corollary \ref{cor1}.
\end{proof}

\bigbreak

\noindent
{\bf Acknowledgment:} The authors thank B. Texier for fruitful discussions.

\bibliographystyle{plain}
\bibliography{bibAlvarezLannes}

\end{document}